\documentclass[12pt]{article}
\usepackage[utf8x]{inputenc}
\usepackage{amssymb, amsmath}
\usepackage[dvips]{graphics}
\usepackage[all]{xy}
\usepackage{bbm}
\usepackage{enumitem}
\usepackage{setspace}
\usepackage[colorlinks=true]{hyperref}
\usepackage{xcolor}
\hypersetup{colorlinks   = true}
\hypersetup{linkcolor=blue}
\begin{document}
\newtheorem{The}{Theorem}[section]
\newtheorem{Lem}[The]{Lemma}
\newtheorem{Prop}[The]{Proposition}
\newtheorem{Cor}[The]{Corollary}
\newtheorem{Rem}[The]{Remark}
\newtheorem{Obs}[The]{Observation}
\newtheorem{SConj}[The]{Standard Conjecture}
\newtheorem{Titre}[The]{\!\!\!\! }
\newtheorem{Conj}[The]{Conjecture}
\newtheorem{Question}[The]{Question}
\newtheorem{Prob}[The]{Problem}
\newtheorem{Def}[The]{Definition}
\newtheorem{Not}[The]{Notation}
\newtheorem{Claim}[The]{Claim}
\newtheorem{Conc}[The]{Conclusion}
\newtheorem{Ex}[The]{Example}
\newtheorem{Fact}[The]{Fact}
\newtheorem{Formula}[The]{Formula}
\newtheorem{Formulae}[The]{Formulae}
\newcommand{\C}{\mathbb{C}}
\newcommand{\R}{\mathbb{R}}
\newcommand{\N}{\mathbb{N}}
\newcommand{\Z}{\mathbb{Z}}
\newcommand{\Q}{\mathbb{Q}}
\newcommand{\Proj}{\mathbb{P}}
\newcommand{\Rc}{\mathcal{R}}
\newcommand{\Oc}{\mathcal{O}}

\begin{center}

{\Large\bf Functionals for the Study of LCK Metrics on Compact Complex Manifolds}

\end{center}

\begin{center}

{\large Dan Popovici and Erfan Soheil}

\end{center}

\vspace{1ex}

\noindent{\small{\bf Abstract.} We propose an approach to the existence problem for locally conformally K\"ahler metrics on compact complex manifolds by introducing and studying a functional that is different according to whether the complex dimension of the manifold is $2$ or higher.}

\vspace{2ex}

\section{Introduction}

Let $X$ be an $n$-dimensional compact complex manifold with $n\geq 2$. In this paper, we propose a variational approach to the existence of locally conformally K\"ahler (lcK) metrics on $X$ by introducing and analysing a functional in each of the cases $n=2$ and $n\geq 3$. This functional, defined on the non-empty set ${\cal H}_X$ of all the Hermitian metrics on $X$, assumes non-negative values and vanishes precisely on the lcK metrics. We compute the first variation of our functional on both surfaces and higher-dimensional manifolds.

\vspace{2ex}

We will identify a Hermitian metric on $X$ with the associated $C^\infty$ positive definite $(1,\,1)$-form $\omega$. The set ${\cal H}_X$ of all these metrics is a non-empty open convex cone in the infinite-dimensional real vector space $C^\infty_{1,\,1}(X,\,\R)$ of all the real-valued smooth $(1,\,1)$-forms on $X$. As is well known, a Hermitian metric $\omega$ is called K\"ahler if $d\omega=0$ and a complex manifold $X$ is said to be K\"ahler if there exists a K\"ahler metric thereon. Meanwhile, the notion of {\bf locally conformally K\"ahler (lcK)} manifold originates with I. Vaisman in [Vai76]. There are several equivalent definitions of lcK manifolds. The one adopted in this paper stipulates that a complex manifold $X$ is lcK if there exists an lcK metric thereon, while a Hermitian metric $\omega$ on $X$ is said to be lcK if there exists a $C^\infty$ $1$-form $\theta$ on $X$ such that $d\theta=0$ and $$d\omega = \omega\wedge\theta.$$ When it exists, the $1$-form $\theta$ is unique and is called the Lee form of $\omega$. For equivalent definitions of lcK manifolds, the reader is referred e.g. to Definitions 3.18 and 3.29 of [OV22].

One of the early results in the theory of lcK manifolds is Vaisman's theorem according to which any lcK metric on a compact K\"ahler manifold is, in fact, globally conformally K\"ahler. This theorem was extended to compact complex spaces with singularities by Preda and Stanciu in [PS22].

The question of when lcK metrics exist on a given compact complex manifold $X$ has been extensively studied. For example, Otiman characterised the existence of such metrics with prescribed Lee form in terms of currents: given a $d$-closed $1$-form $\theta$ on $X$ and considering the associated twisted operator $d_{\theta}=d +\theta \wedge\cdot$, Theorem 2.1 in [Oti14] stipulates that $X$ admits an lcK metric whose Lee form is $\theta$ if and only if there are no non-trivial positive $(1, 1)$-currents on $X$ that are $(1, 1)$-components of $d_{\theta}$-boundaries.

On the other hand, Istrati investigated the relation between the existence of special lcK metrics on a compact complex manifold and the group of biholomorphisms of the manifold. Specifically, according to Theorem 0.2 in [Ist19], a compact lcK manifold $X$ admits a Vaisman metric if the group of biholomorphisms of $X$ contains a torus $\mathbb{T}$ that is not purely real. A compact torus $\mathbb{T}$ of biholomorphisms of a compact complex manifold $(X,\,J)$ is said to be {\it purely real} (in the sense of (1) of Definition 0.1. in [Ist19]) if its Lie algebra $\mathfrak{t}$ satisfies the condition $\mathfrak{t}\cap J\mathfrak{t}= 0$, where $J$ is the complex structure of $X$. Recall that an lcK metric $\omega$ is said to be a {\it Vaisman metric} if $\nabla^\omega\theta=0$, where $\theta$ is the Lee form of $\omega$ and $\nabla^\omega$ is the Levi-Civita connection determined by $\omega$.

\vspace{2ex}

The approach we propose in this paper to the issue of the existence of lcK metrics on a compact complex $n$-dimensional manifold $X$ is analytic. Given an arbitrary Hermitian metric $\omega$ on $X$, the Lefschetz decomposition \begin{eqnarray*}d\omega = (d\omega)_{prim} + \omega\wedge\theta_\omega\end{eqnarray*} of $d\omega$ into a uniquely determined $\omega$-primitive part and a part divisible by $\omega$ with a uniquely determined quotient $1$-form $\theta_\omega$ (the Lee form of $\omega$) gives rise to the following dichotomy (cf. Lemma \ref{Lem:general-obstruction}):

  \vspace{1ex}

  (i)\, either $n=2$, in which case $(d\omega)_{prim} = 0$ but the Lee form $\theta_\omega$ need not be $d$-closed, so the lcK condition on $\omega$ is equivalent to $d\theta_\omega = 0$. This turns out to be equivalent to $\partial\theta_\omega^{1,\,0} = 0$. Therefore, we define our functional $L:{\cal H}_X\longrightarrow[0,\,+\infty)$ in this case to be $$L(\omega) = ||\partial\theta_\omega^{1,\,0}||_\omega^2,$$ namely its value at every Hermitian metric $\omega$ on $X$ is defined to be the squared $L^2_\omega$-norm of $\partial\theta_\omega^{1,\,0}$.

\vspace{1ex}

  (ii)\, or $n\geq 3$, in which case the lcK condition on $\omega$ is equivalent to the vanishing condition $(d\omega)_{prim} = 0$. This is further equivalent to the vanishing of either $(\partial\omega)_{prim}$ or $(\bar\partial\omega)_{prim}$. We, therefore, define our functional $L:{\cal H}_X\longrightarrow[0,\,+\infty)$ in this case to be $$L(\omega) = ||(\bar\partial\omega)_{prim}||_\omega^2,$$ namely its value at every Hermitian metric $\omega$ on $X$ is defined to be the squared $L^2_\omega$-norm of the $\omega$-primitive part of the $(1,\,2)$-form $\bar\partial\omega$.

  \vspace{2ex}

  In [Vai90], Vaisman considered two functionals in complex dimension $2$, one of which is a variant of our functional $L$ for complex surfaces. It actually coincides with the auxiliary functional ${\cal L}$ that we consider in the proof of Lemma \ref{Lem:functional_justification}. Vaisman computes both the first and second variation of his functionals, but he employs a different method to ours. As far as we are aware, there is yet no analogue to our functional in complex dimension $>2$.

 \vspace{2ex}

The main results of the paper are the computations of the first variation of our functional $L$ in each of the cases $n=2$ (cf. Theorem \ref{The:1st-variation_surfaces}) and $n\geq 3$ (cf. Theorem \ref{The:derivative_L_deg>2}). From these, we deduce the Euler-Lagrange equations (cf. Corollaries \ref{Cor:1st-variation_surfaces_cor_E-L_equation} and \ref{Cor:E-L_equation_dim>2}). The equation for $n=2$ has an analogue in [Vai90], while the one for $n\geq 3$ seems entirely new.

While the functional $L$ is scaling-invariant when $n=2$, this fails to be the case when $n\geq 3$. In this latter case, we obtain two proofs -- one as a corollary of the formula for the first variation of our functional (cf. Proposition \ref{Prop:lcK-critical}), the other as a direct consequence of the behaviour of our functional in the scaling direction (cf. Proposition \ref{Prop:lcK-critical_bis}) -- for the equivalence:

\vspace{1ex}

\hspace{6ex}  $\omega$ is a {\bf critical point} for the functional $L$ if and only if $\omega$ is {\bf lcK}

\vspace{1ex}

Still in the case $n\geq 3$, we introduce in Definition \ref{Def:normalised-functional_lcK} a normalised version $\widetilde{L}_\rho$ of the functional $L$ depending on an arbitrary background Hermitian metric $\rho$. The first variation of $\widetilde{L}_\rho$ is then deduced in Proposition \ref{Prop:1st-variation-normalised-functional_lcK} from the analogous computation for $L$ obtained in Theorem \ref{The:derivative_L_deg>2}. One motivation for the normalisation we propose in terms of a (possibly balanced and possibly moving) metric $\rho$ stems from the conjecture predicting that the simultaneous existence of a balanced metric and of an lcK metric on a compact complex manifold ought to imply the existence of a K\"ahler metric. We hope to be able to develop this line of thought in future work.

At the end of $\S.$\ref{section:normalised-functionals}, we use our scaling-invariant functionals $L$ (in the case of compact complex surfaces) and $\widetilde{L}_\rho$ (in the case of higher-dimensional compact complex manifolds) to produce positive $(1,\,1)$-currents whose failure to be either $C^\infty$ forms or strictly positive provides possible obstructions to the existence of lcK metrics.

\vspace{2ex}

\noindent {\bf Acknowledgments.} This work is part of the second-named author's thesis under the supervision of the first-named author. The former wishes to thank the latter for constant support. Both authors are very grateful to the referee for their careful reading of the text, their helpful remarks and suggestions, as well as for pointing out the reference [Vai90] the authors were unfortunately unaware of when writing the first version of this paper.

\section{Preliminaries}\label{section:standard}

In this section, we recast some standard material in the language of primitive forms and make a few observations that will be used in the next sections. 

Let $X$ be a complex manifold with $\mbox{dim}_\C X = n$. We will denote by:

\vspace{1ex}

(i)\, $C^\infty_k(X,\,\C)$, resp. $C^\infty_{p,\,q}(X,\,\C)$, the space of $C^\infty$ differential forms of degree $k$, resp. of bidegree $(p,\,q)$ on $X$. When these forms $\alpha$ are real (in the sense that $\overline\alpha = \alpha$), the corresponding spaces will be denoted by $C^\infty_k(X,\,\R)$, resp. $C^\infty_{p,\,q}(X,\,\R)$.

\vspace{1ex}

(ii)\, $\Lambda^kT^\star X$, resp. $\Lambda^{p,\,q}T^\star X$, the vector bundle of differential forms of degree $k$, resp. of bidegree $(p,\,q)$, as well as the spaces of such forms considered in a pointwise way.

\vspace{1ex}

For any $(1,\,1)$-form $\rho\geq 0$, we will also use the following notation: \begin{eqnarray*}\rho_{k}:= \frac{\rho^{k}}{k!}, \hspace{5ex} 1\leq k \leq n.\end{eqnarray*} When $\rho=\omega$ is $C^\infty$ and positive definite (i.e. $\omega$ is a Hermitian metric on $X$), it can immediately be checked that \begin{eqnarray*}d\omega_{k} = \omega_{k-1}\wedge d\omega \hspace{3ex}\mbox{and}\hspace{3ex} \star_{\omega}\omega_{k}=\omega_{n-k} \end{eqnarray*} for all $1\leq k\leq n$, where $\star = \star_\omega$ is the Hodge star operator induced by $\omega$.

\vspace{2ex}  

Recall the following standard

\begin{Def} A $C^\infty$ positive definite $(1,\,1)$-form (i.e. a Hermitian metric) $\omega$ on a complex manifold $X$ is said to be {\bf locally conformally K\"ahler (lcK)} if 

$$d\omega = \omega\wedge\theta  \hspace{3ex}  \mbox{for some} \hspace{1ex} C^\infty \hspace{1ex} 1\mbox{-form} \hspace{1ex} \theta  \hspace{1ex}  \mbox{satisfying} \hspace{1ex} d\theta =0.$$

\noindent The $1$-form $\theta$ is uniquely determined, is real and is called the {\bf Lee form} of $\omega$.

\end{Def}

\noindent {\it (A)\, Hermitian-geometric preliminaries}

\vspace{1ex}

$\bullet$ Recall that for any $k\leq n$ and any Hermitian metric $\omega$ on $X$, the multiplication map 

$$L_\omega^l=\omega^l\wedge\cdot : \Lambda^kT^\star X\longrightarrow\Lambda^{k+2l}T^\star X$$

\noindent defined at every point of $X$ is an {\bf isomorphism} if $l=n-k$, is {\bf injective} (but in general not surjective) for every $l< n-k$ and is {\bf surjective} (but in general not injective) for every $l>n-k$. A $k$-form is said to be $\omega$-{\bf primitive} if it lies in the kernel of the multiplication map $L^{n-k+1}_\omega$. Equivalently, the $\omega$-primitive $k$-forms are precisely those that lie in the kernel of $\Lambda_\omega : \Lambda^k T^\star X\longrightarrow\Lambda^{k-2}T^\star X$, the adjoint w.r.t. the pointwise inner product $\langle\cdot,\cdot\rangle_\omega$ (hence also w.r.t. the $L^2$-inner product $\langle\langle\cdot,\cdot\rangle\rangle_\omega$) of the Lefschetz operator $L_\omega=\omega\wedge\cdot : \Lambda^kT^\star X\longrightarrow\Lambda^{k+2}T^\star X$.

\vspace{2ex}

$\bullet$ Also recall that for every $k\leq n$, every $k$-form $\alpha$ admits a unique $\langle \,\,,\,\,\rangle_\omega$-orthogonal pointwise splitting (called the {\bf Lefschetz decomposition}): \begin{eqnarray}\label{eqn:Lefschetz-decomp}\alpha = \alpha_{prim} + \omega\wedge\beta_{prim}^{(1)} + \omega^2\wedge\beta_{prim}^{(2)} + \dots  + \omega^r\wedge\beta_{prim}^{(r)},\end{eqnarray} where $r$ is the largest non-negative integer such that $2r\leq k$, $\alpha_{prim}$, $\beta_{prim}^{(1)}, \dots , \beta_{prim}^{(r)}$ are $\omega$-primitive forms of respective degrees $k$, $k-2, \dots , k-2r\geq 0$, and $\langle \,\,,\,\,\rangle_\omega$ is the pointwise inner product defined by $\omega$. We will call $\alpha_{prim}$ the {\it primitive part} of $\alpha$.

\vspace{2ex}

$\bullet$ The following general formula (cf. e.g. [Voi02, Proposition 6.29, p. 150]) that holds for any {\it primitive} form $v$ of arbitrary bidegree $(p, \, q)$ on any complex $n$-dimensional manifold will be of great use: \begin{eqnarray}\label{eqn:prim-form-star-formula-gen}\star\, v = (-1)^{k(k+1)/2}\, i^{p-q}\,\omega_{n-p-q}\wedge v, \hspace{2ex} \mbox{where}\,\, k:=p+q,\end{eqnarray}

\vspace{2ex}

 $\bullet$ We will often use the standard notation $[A,\,B]:=AB-(-1)^{ab}BA$, where $A$ and $B$ are arbitrary linear operators of degrees $a$, resp. $b$, acting on the differential forms of $X$. The following formula (see e.g. [Dem97, VI, $\S5.2$, Corollary 5.9]) will come in handy several times: \begin{eqnarray}\label{eqn:L-Lambda_commutation}[\Lambda_\omega,\,L_\omega] = (n-k)\,\mbox{Id}\end{eqnarray} when acting on $k$-forms on $X$.

\vspace{2ex}

 $\bullet$ Finally, recall the Hermitian commutation relation: \begin{equation}\label{eqn:Hermitian-commutation}i[\Lambda_\omega,\,\partial] = -(\bar\partial_\omega^\star + \bar\tau_\omega^\star)\end{equation} proved in [Dem84], where $\tau_\omega:=[\Lambda_\omega,\,\partial\omega\wedge\cdot]$ is the torsion operator of order $0$ and bidegree $(1,\,0)$. This definition of $\tau_\omega$ yields

\vspace{1ex}

\hspace{25ex} $\bar\tau_\omega^\star\omega = [(\bar\partial\omega\wedge\cdot)^\star, L_\omega](\omega) = (\bar\partial\omega\wedge\cdot)^\star(\omega^2).$

\vspace{2ex}

$\bullet$ On the other hand, if $\alpha^{1,\,0}$ is any $(1,\,0)$-form on $X$, let $\xi_{\alpha^{1,\,0}}$ be the $(1,\,0)$-vector field defined by the requirement $\bar\xi_{\alpha^{1,\,0}}\lrcorner\omega = \alpha^{1,\,0}$. If we set $\alpha^{0,\,1}:=\overline{\alpha^{1,\,0}}$, we have $\bar\xi_{\alpha^{1,\,0}} = \xi_{\alpha^{0,\,1}}$, where $\xi_{\alpha^{0,\,1}}$ is the $(0,\,1)$-vector field defined by the requirement $\bar\xi_{\alpha^{0,\,1}}\lrcorner\omega = \alpha^{0,\,1}$.

It is easily checked in local coordinates chosen about a given point $x$ such that the metric $\omega$ is defined by the identity matrix at $x$, that the adjoint w.r.t. $\langle\,\,,\,\,\rangle_\omega$ of the contraction operator by $\bar\xi_{\alpha^{1,\,0}} = \xi_{\alpha^{0,\,1}}$ is given by the formula \begin{equation}\label{eqn:xi_bar_contracting_theta01_pre_first}(\bar\xi_{\alpha^{1,\,0}}\lrcorner\cdot)^\star = (\xi_{\alpha^{0,\,1}}\lrcorner\cdot)^\star  = -i\alpha^{0,\,1}\wedge\cdot, \hspace{3ex} \mbox{or equivalently} \hspace{3ex} -i\bar\xi_{\alpha^{1,\,0}}\lrcorner\cdot = -i\xi_{\alpha^{0,\,1}}\lrcorner\cdot = (\alpha^{0,\,1}\wedge\cdot)^\star.\end{equation}

\noindent Taking conjugates, we get: \begin{equation}\label{eqn:xi_bar_contracting_theta01_pre}(\xi_{\alpha^{1,\,0}}\lrcorner\cdot)^\star = i\alpha^{1,\,0}\wedge\cdot, \hspace{3ex} \mbox{or equivalently} \hspace{3ex} i\xi_{\alpha^{1,\,0}}\lrcorner\cdot = (\alpha^{1,\,0}\wedge\cdot)^\star.\end{equation}

Explicitly, if $\alpha^{0,\,1} = \sum\limits_k\bar{a}_kd\bar{z}_k$ on a neighbourhood of $x$, then $ -i\bar\xi_\alpha\lrcorner\cdot = (\alpha^{0,\,1}\wedge\cdot)^\star =  \sum\limits_k a_k\,\frac{\partial}{\partial\bar{z}_k}\lrcorner\cdot$ at $x$. Hence, $-i\bar\xi_\alpha\lrcorner\alpha^{0,\,1} = \sum\limits_k|a_k|^2 = |\alpha^{0,\,1}|^2_\omega$ at $x$. We have just got the pointwise formula: \begin{equation}\label{eqn:xi_bar_contracting_theta01}-i\bar\xi_\alpha\lrcorner\alpha^{0,\,1} = |\alpha^{0,\,1}|^2_\omega = |\alpha^{1,\,0}|^2_\omega\end{equation}

  \noindent at every point of $X$.

\vspace{2ex}

 \noindent {\it (B)\, lcK-geometric preliminaries}

\vspace{1ex}

 Now, suppose that $d\omega=\omega\wedge\theta_\omega$ for some (necessarily real) $1$-form $\theta_\omega$. Then, $\bar\partial\omega=\omega\wedge\theta_\omega^{0,\,1}$, so $(\bar\partial\omega\wedge\cdot)^\star = -i\Lambda_\omega(\bar\xi_\theta\lrcorner\cdot)$, where $\bar\xi_\theta:=\bar\xi_\alpha$ with $\alpha^{1,\,0}=\theta_\omega^{1,\,0}$. The above formula for $\bar\tau_\omega^\star\omega$ translates to

$$\bar\tau_\omega^\star\omega =-i\Lambda_\omega(\bar\xi_\theta\lrcorner\omega^2) = -2i\Lambda_\omega(\omega\wedge(\bar\xi_\theta\lrcorner\omega)) = -2i[\Lambda_\omega,\,L_\omega](\bar\xi_\theta\lrcorner\omega) = -2i(n-1)\theta_\omega^{1,\,0}$$

 The conclusion of this discussion is that, when $d\omega=\omega\wedge\theta_\omega$, formula (\ref{eqn:Lee-form_formula_higher}) (which will be proved as part of Lemma \ref{Lem:general-obstruction} below) translates to $$\theta_\omega^{1,\,0} = \frac{1}{n-1}\,\Lambda_\omega(\partial\omega) = \frac{1}{n-1}\,[\Lambda_\omega,\,\partial](\omega) = \frac{1}{n-1}\,i\bar\partial^\star_\omega \omega + \frac{1}{n-1}\,i\bar\tau_\omega^\star\omega = \frac{1}{n-1}\,i\bar\partial^\star_\omega \omega + 2\theta_\omega^{1,\,0},$$ \noindent which amounts to $\theta_\omega^{1,\,0} = -\frac{1}{n-1}\,i\bar\partial^\star_\omega \omega.$ This proves (\ref{eqn:Lee-form_formula_1-0_higher}) for an arbitrary $n$, hence also (\ref{eqn:Lee-form_formula_1-0}) when $n=2$, if the other statements in Lemma \ref{Lem:general-obstruction} have been proved.

 \vspace{2ex}

The obstruction to a given Hermitian metric $\omega$ being lcK depends on whether $n=2$ or $n\geq 3$.

\begin{Lem}\label{Lem:general-obstruction} Let $X$ be a complex manifold with $\mbox{dim}_\C X=n$.

\vspace{1ex}

$(i)$\, If $n=2$, for any Hermitian metric $\omega$ there exists a unique, possibly non-closed, $C^\infty$ $1$-form $\theta = \theta_\omega$ such that $d\omega = \omega\wedge\theta$. Therefore, $\omega$ is {\bf lcK} if and only if $\theta_\omega$ is {\bf $d$-closed}. 

Moreover, for any Hermitian metric $\omega$, the $2$-form $d\theta_\omega$ is $\omega$-{\bf primitive}, i.e. $\Lambda_\omega(d\theta_\omega)=0$, or equivalently, $\omega\wedge d\theta_\omega=0$, while the Lee form is real and is explicitly given by the formula: \begin{equation}\label{eqn:Lee-form_formula}\theta_\omega = \Lambda_\omega(d\omega).\end{equation}

\noindent Alternatively, if $\theta_\omega = \theta_\omega^{1,\,0} + \theta_\omega^{0,\,1}$ is the splitting of $\theta_\omega$ into components of pure types, we have \begin{equation}\label{eqn:Lee-form_formula_1-0}\theta_\omega^{1,\,0} = \Lambda_\omega(\partial\omega) = -i\bar\partial^\star\omega\end{equation}

\noindent and the analogous formulae for $\theta_\omega^{0,\,1} = \overline{\theta_\omega^{1,\,0}}$ obtained by taking conjugates. 

\vspace{1ex}

\noindent $(ii)$\, If $n\geq 3$, for any Hermitian metric $\omega$ there exists a unique $\omega$-primitive $C^\infty$ $3$-form $(d\omega)_{prim}$ and a unique $C^\infty$ $1$-form $\theta = \theta_\omega$ such that $d\omega = (d\omega)_{prim} + \omega\wedge\theta$. The Lee form is real and is explicitly given by the formula \begin{equation}\label{eqn:Lee-form_formula_higher}\theta_\omega = \frac{1}{n-1}\,\Lambda_\omega(d\omega).\end{equation}

\noindent Moreover, $\omega$ is {\bf lcK} if and only if {\bf $(d\omega)_{prim} = 0$}.

 If $\omega$ is lcK, then \begin{equation}\label{eqn:Lee-form_formula_1-0_higher}\theta_\omega^{1,\,0} = \frac{1}{n-1}\,\Lambda_\omega(\partial\omega) = -\frac{i}{n-1}\,\bar\partial^\star\omega\end{equation}

\noindent and the analogous formulae obtained by taking conjugates hold for $\theta_\omega^{0,\,1} = \overline{\theta_\omega^{1,\,0}}$.

\end{Lem}

\noindent {\it Proof.} $(i)$\, When $n=2$, the map $\omega\wedge\cdot : \Lambda^1T^\star X\longrightarrow\Lambda^3 T^\star X$ is an isomorphism at every point of $X$. In particular, the $3$-form $d\omega$ is the image of a unique $1$-form $\theta$ under this map. 

 To see that $d\theta$ is primitive, we apply $d$ to the identity $d\omega = \omega\wedge\theta$ to get

$$0=d^2\omega = d\omega\wedge\theta + \omega\wedge d\theta.$$

\noindent Meanwhile, multiplying the same identity by $\theta$, we get $d\omega\wedge\theta = \omega\wedge\theta\wedge\theta = 0$ since $\theta\wedge\theta=0$ due to the degree of $\theta$ being $1$. Therefore, $\omega\wedge d\theta=0$, which means that the $2$-form $d\theta$ is $\omega$-primitive.

 To prove formula (\ref{eqn:Lee-form_formula}), we apply $\Lambda_\omega$ to the identity $d\omega = \omega\wedge\theta$ to get

$$\Lambda_\omega(d\omega) = [\Lambda_\omega,\,L_\omega](\theta) = \theta,$$

\noindent where we used the identities $\Lambda_\omega(\theta)=0$ (for bidegree reasons) and (\ref{eqn:L-Lambda_commutation}) (with $k=1$ and $n=2$).

\vspace{1ex}

$(ii)$\, The splitting $d\omega = (d\omega)_{prim} + \omega\wedge\theta$ is the Lefschetz decomposition of $d\omega$ w.r.t. the metric $\omega$. Applying $\Lambda_\omega$, we get $\Lambda_\omega(d\omega) = [\Lambda_\omega,\,L_\omega](\theta) = (n-1)\,\theta$ (having applied (\ref{eqn:L-Lambda_commutation}) with $k=1$ to get the latter identity), which proves (\ref{eqn:Lee-form_formula_higher}).

 The implication ``$\omega$ lcK $\implies (d\omega)_{prim}=0$`` follows at once from the definitions. To prove the reverse implication, suppose that $(d\omega)_{prim}=0$. We have to show that $\theta$ is $d$-closed. The assumption means that $d\omega = \omega\wedge\theta$, so $d\omega\wedge\theta = \omega\wedge\theta\wedge\theta=0$ and $0=d^2\omega = d\omega\wedge\theta + \omega\wedge d\theta$. Consequently, $\omega\wedge d\theta =0$. Now, the multiplication of $k$-forms by $\omega^l$ is injective whenever $l\leq n-k$. When $n\geq 3$, if we choose $l=1$ and $k=2$ we get that the multiplication of $2$-forms by $\omega$ is injective. Hence, the identity $\omega\wedge d\theta =0$ implies $d\theta =0$, so $\omega$ is lcK.   \hfill $\Box$

\vspace{3ex}

Another standard observation is that the Lefschetz decomposition transforms nicely, hence the lcK property is preserved, under conformal rescaling.

\begin{Lem}\label{Lem:conformal-rescaling} Let $\omega$ be an arbitrary Hermitian metric and let $f$ be any smooth real-valued function on a compact complex $n$-dimensional manifold $X$. If $d\omega = (d\omega)_{prim} + \omega\wedge\theta_\omega$ is the Lefschetz decomposition of $d\omega$ w.r.t. the metric $\omega$ (with the understanding that $(d\omega)_{prim} = 0$ when $n=2$), then

\begin{equation}\label{eqn:rescaled-Lefschetz-decomp}d(e^f\omega) = e^f(d\omega)_{prim} + e^f\omega\wedge(\theta_\omega + df)\end{equation}

\noindent is the Lefschetz decomposition of $d(e^f\omega)$ w.r.t. the metric $\widetilde\omega:=e^f\omega$.

 Consequently, $\omega$ is lcK if and only if any conformal rescaling $e^f\omega$ of $\omega$ is lcK, while the Lee form transforms as $\theta_{e^f\omega} = \theta_\omega + df$. In particular, when the lcK metric $\omega$ varies in a fixed conformal class, the Lee form $\theta_\omega$ varies in a fixed De Rham $1$-class $\{\theta_\omega\}_{DR}\in H^1(X,\,\R)$ called the {\bf Lee De Rham class} associated with the given conformal class. Moreover, the map $\omega\mapsto\theta_\omega$ defines a bijection from the set of lcK metrics in a given conformal class to the set of elements of the corresponding Lee De Rham $1$-class.   

\end{Lem}

\noindent {\it Proof.} Differentiating, we get $ d(e^f\omega) = e^f d\omega + e^f\omega\wedge df = e^f(d\omega)_{prim} + e^f\omega\wedge(\theta_\omega + df)$. Meanwhile, it can immediately be checked that 

$$\Lambda_{e^f\omega} = e^{-f}\Lambda_\omega,$$

\noindent so $\ker\Lambda_{e^f\omega} = \ker\Lambda_\omega$. Thus, the $\omega$-primitive forms coincide with the $\widetilde\omega$-primitive forms. Since $\Lambda_{\widetilde\omega}$ commutes with the multiplication by any real-valued function, $ e^f(d\omega)_{prim}$ is $\widetilde\omega$-primitive, so (\ref{eqn:rescaled-Lefschetz-decomp}) is the Lefschetz decompostion of $d\widetilde\omega$ w.r.t. $\widetilde\omega$. \hfill $\Box$

\vspace{3ex}

 When $X$ is compact, we know from [Gau77] that every Hermitian metric $\omega$ on $X$ admits a (unique up to a positive multiplicative constant) conformal rescaling $\widetilde\omega:=e^f\omega$ that is a {\bf Gauduchon metric}. These metrics are defined (cf. [Gau77]) by the requirement that $\partial\bar\partial\widetilde\omega^{n-1}=0$, where $n$ is the complex dimension of $X$. This fact, combined with Lemma \ref{Lem:conformal-rescaling}, shows that no loss of generality is incurred in the study of the existence of lcK metrics on compact complex manifolds if we confine ourselves to Gauduchon metrics.

 We end this review of known material with the following characterisation (cf. [AD15, Lemma 2.5]) of Gauduchon metrics on surfaces in terms of their Lee forms. It appears that, in any dimension, a metric $\omega$ is Gauduchon if and only if $d^\star_\omega\theta_\omega=0$\footnote{The authors are grateful to the referee for pointing this fact out to them.}, but we confine ourselves to the $2$-dimensional case.

\begin{Lem}\label{Lem:Gauduchon_Lee-form_characterisation} Let $\omega$ be a Hermitian metric on a complex {\bf surface} $X$. The following equivalence holds:

$$\partial\bar\partial\omega = 0 \hspace{2ex} \mbox{(i.e.}\hspace{1ex} \omega\hspace{1ex} \mbox{is a Gauduchon metric)} \hspace{2ex} \iff \hspace{2ex} \bar\partial^\star_\omega\theta_\omega^{0,\,1} =0,$$

\noindent where $\theta_\omega^{0,\,1}$ is the component of type $(0,\,1)$ of the Lee form $\theta_\omega$ of $\omega$. 

\end{Lem}

\noindent {\it Proof.} We give a proof different from the one in [AD15] by making use of the Hermitian commutation relations. By applying $\partial$ to the identity $\bar\partial\omega = \omega\wedge\theta_\omega^{0,\,1}$ and using the identity $\partial\omega = \omega\wedge\theta_\omega^{1,\,0}$, we get

$$\partial\bar\partial\omega = \partial\omega\wedge\theta_\omega^{0,\,1} + \omega\wedge\partial\theta_\omega^{0,\,1} = \omega\wedge(\theta_\omega^{1,\,0}\wedge\theta_\omega^{0,\,1} + \partial\theta_\omega^{0,\,1}).$$

\noindent Taking $\Lambda_\omega$, we get $$\Lambda_\omega(\partial\bar\partial\omega) = [\Lambda_\omega,\,L_\omega](\theta_\omega^{1,\,0}\wedge\theta_\omega^{0,\,1} + \partial\theta_\omega^{0,\,1}) + \omega\wedge\Lambda_\omega(\theta_\omega^{1,\,0}\wedge\theta_\omega^{0,\,1} + \partial\theta_\omega^{0,\,1}) = \Lambda_\omega(\theta_\omega^{1,\,0}\wedge\theta_\omega^{0,\,1} + \partial\theta_\omega^{0,\,1})\,\omega,$$

\noindent where the second identity follows from $[\Lambda_\omega,\,L_\omega] = -(2-2)\,\mbox{Id} = 0$ on $2$-forms on complex surfaces. Now, $\Lambda_\omega(\theta_\omega^{1,\,0}\wedge\theta_\omega^{0,\,1} + \partial\theta_\omega^{0,\,1})$ is a function, so from the above identities we get the equivalences \begin{eqnarray}\nonumber \Lambda_\omega(\partial\bar\partial\omega) = 0 & \iff & \Lambda_\omega(\theta_\omega^{1,\,0}\wedge\theta_\omega^{0,\,1} + \partial\theta_\omega^{0,\,1})=0 \iff \theta_\omega^{1,\,0}\wedge\theta_\omega^{0,\,1} + \partial\theta_\omega^{0,\,1} \hspace{1ex} \mbox{is $\omega$-primitive} \\
\nonumber & \iff & \omega\wedge(\theta_\omega^{1,\,0}\wedge\theta_\omega^{0,\,1} + \partial\theta_\omega^{0,\,1}) = 0 \iff \partial\bar\partial\omega=0.\end{eqnarray}

 We remember the equivalence $\partial\bar\partial\omega=0 \iff \Lambda_\omega(\theta_\omega^{1,\,0}\wedge\theta_\omega^{0,\,1}) + \Lambda_\omega(\partial\theta_\omega^{0,\,1})=0$. Since $\Lambda_\omega(i\theta_\omega^{1,\,0}\wedge\theta_\omega^{0,\,1}) = |\theta_\omega^{1,\,0}|_\omega^2$ (immediate verification) and $\Lambda_\omega\theta_\omega^{0,\,1} = 0$ (for bidegree reasons), we get the equivalence:

$$\partial\bar\partial\omega=0 \iff |\theta_\omega^{1,\,0}|_\omega^2 + i[\Lambda_\omega,\,\partial]\,\theta_\omega^{0,\,1}=0.$$

 The Hermitian commutation relation $i[\Lambda_\omega,\,\partial] = -(\bar\partial_\omega^\star + \bar\tau_\omega^\star)$ (cf. (\ref{eqn:Hermitian-commutation}), see [Dem84]) transforms the last equivalence into

\begin{equation}\label{eqn:equivalence-before-formula}\partial\bar\partial\omega=0 \iff |\theta_\omega^{1,\,0}|_\omega^2 - (\bar\partial_\omega^\star\theta_\omega^{0,\,1} + \bar\tau_\omega^\star\theta_\omega^{0,\,1}) = 0.\end{equation}

 On the other hand, $\bar\tau_\omega^\star = [(\bar\partial\omega\wedge\cdot)^\star,\,\omega\wedge\cdot]$. From this we get

\begin{Formula}\label{Formula:tau_bar_star-theta01} For any Hermitian metric $\omega$ on a complex surface, we have

$$\bar\tau_\omega^\star\theta_\omega^{0,\,1} = |\theta_\omega^{0,\,1}|^2_\omega.$$

\end{Formula}

\noindent {\it Proof of Formula \ref{Formula:tau_bar_star-theta01}.} Since $(\bar\partial\omega\wedge\cdot)^\star\theta_\omega^{0,\,1}=0$ for bidegree reasons, we get $\bar\tau_\omega^\star\theta_\omega^{0,\,1} = (\bar\partial\omega\wedge\cdot)^\star(\omega\wedge\theta_\omega^{0,\,1})$. Since $\bar\partial\omega=\omega\wedge\theta_\omega^{0,\,1}$, we have $(\bar\partial\omega\wedge\cdot)^\star = -i\Lambda_\omega(\bar\xi_\theta\lrcorner\cdot)$ (see (\ref{eqn:xi_bar_contracting_theta01}) and the discussion there below), where $\bar\xi_\theta$ is the $(0,\,1)$-vector field defined by the requirement $\bar\xi_\theta\lrcorner\omega=\theta_\omega^{1,\,0}$. Hence

$$\bar\tau_\omega^\star\theta_\omega^{0,\,1} = -i\Lambda_\omega(\theta_\omega^{1,\,0}\wedge\theta_\omega^{0,\,1}) -i\Lambda_\omega[\omega\wedge(\bar\xi_\theta\lrcorner\theta_\omega^{0,\,1})].$$

\noindent Since $-i\bar\xi_\theta\lrcorner\theta_\omega^{0,\,1} = |\theta_\omega^{0,\,1}|^2_\omega$ (cf. (\ref{eqn:xi_bar_contracting_theta01})), we infer that $$\bar\tau_\omega^\star\theta_\omega^{0,\,1} = -\Lambda_\omega(i\theta_\omega^{1,\,0}\wedge\theta_\omega^{0,\,1}) + 2\,|\theta_\omega^{0,\,1}|^2_\omega,$$

\noindent since $\Lambda_\omega(\omega)=n=2$. Meanwhile, $\theta_\omega^{1,\,0}=\overline{\theta_\omega^{0,\,1}}$, so we get $\Lambda_\omega(i\theta_\omega^{1,\,0}\wedge\theta_\omega^{0,\,1})=|\theta_\omega^{1,\,0}|^2_\omega = |\theta_\omega^{0,\,1}|^2_\omega$ (immediate verification in local coordinates). Formula \ref{Formula:tau_bar_star-theta01} is now proved. \hfill $\Box$

\vspace{3ex}

\noindent {\it End of proof of Lemma \ref{Lem:Gauduchon_Lee-form_characterisation}.} Formula \ref{Formula:tau_bar_star-theta01} transforms equivalence (\ref{eqn:equivalence-before-formula}) into

$$\partial\bar\partial\omega=0 \iff (|\theta_\omega^{1,\,0}|_\omega^2 - |\theta_\omega^{0,\,1}|^2_\omega) - \bar\partial_\omega^\star\theta_\omega^{0,\,1} =0 \iff \bar\partial_\omega^\star\theta_\omega^{0,\,1} =0$$

\noindent and we are done  \hfill $\Box$

\section{An enerygy functional for the study of lcK metrics}\label{section:functional}

%Let $X$ be a compact complex manifold with $\mbox{dim}_\C X=n$. Thanks to Lemma \ref{Lem:conformal-rescaling}, we may restrict attention to Gauduchon metrics in our quest for lcK metrics on $X$. Now, any Gauduchon metric $\omega$ on $X$ represents an Aeppli cohomology class $[\omega^{n-1}]_A\in H^{n-1,\,n-1}(X,\,\R)$, while any positive definite, smooth $(n-1,\,n-1)$-form $\Omega>0$ admits a unique $(n-1)^{st}$ root (i.e. there exists a unique positive definite $(1,\,1)$-form $\omega>0$ such that $\omega^{n-1}=\Omega$) -- see e.g. [Mic83]. So, giving $\omega>0$ in bidegree $(1,\,1)$ is equivalent to giving $\Omega$ in bidegree $(n-1,\,n-1)$. 

%Therefore, it seems reasonable to search for an lcK metric in the Aeppli cohomology class of any Gauduchon metric on $X$. Given a {\bf Gauduchon metric} $\omega_0$, 
In what follows, we will restrict attention to the set $${\cal H}_{X}:=\{\omega\in C^\infty_{1,\,1}(X,\,\R)\,\mid\, \omega>0 \}$$ of all Hermitian metrics on $X$. This is a non-empty open cone in the infinite-dimensional vector space $C^{\infty} _{1,\,1} (X,\, \R)$ of all smooth real $(1,\,1)$-forms on $X$. It will be called the {\bf Hermitian cone} of $X$.

%\noindent be the set of all Hermitian metrics $\omega$ for which $\omega^{n-1}$ lies in the same Aeppli cohomology class as $\omega_0^{n-1}$.

Building on Lemma \ref{Lem:general-obstruction}, we introduce the following energy functional. By $||\,\,\,||_\omega$, respectively $|\,\,\,|_\omega$, we mean the $L^2$-norm, respectively the pointwise norm, defined by $\omega$.

\begin{Def}\label{Def:functional_lcK} Let $X$ be a compact complex manifold with $\mbox{dim}_\C X=n$. 

 $(i)$\, If $n=2$, let $L:{\cal H}_X \longrightarrow[0,\,+\infty)$ be defined by $$L(\omega):=\int\limits_X\partial\theta_\omega^{1,\,0}\wedge\bar\partial\theta_\omega^{0,\,1} = ||\partial\theta_\omega^{1,\,0}||^2_\omega,$$ where $\theta_\omega$ is the Lee form of $\omega$.

\vspace{1ex}

$(ii)$\, If $n\geq 3$, let $L:{\cal H}_X \longrightarrow[0,\,+\infty)$ be defined by $$L(\omega):=\int\limits_X i(\bar\partial\omega)_{prim}\wedge\overline{(\bar\partial\omega)_{prim}}\wedge\omega_{n-3} = ||(\bar\partial\omega)_{prim}||^2_\omega,$$ where $(\bar\partial\omega)_{prim}$ is the $\omega$-primitive part of $\bar\partial\omega$ in its Lefschetz decomposition (\ref{eqn:Lefschetz-decomp}).

\end{Def}

This definition is justified by the following observation.

\begin{Lem}\label{Lem:functional_justification} In the setup of Definition \ref {Def:functional_lcK}, for every metric $\omega\in{\cal H}_X$ the following equivalence holds: $$\omega \hspace{2ex} \mbox{is an lcK metric} \iff L(\omega)=0.$$

\end{Lem}

\noindent {\it Proof.} $\bullet$ In the case $n=2$, we know from $(i)$ of Lemma \ref{Lem:general-obstruction} that $\omega$ is lcK if and only if $d\theta_\omega=0$. This condition is equivalent to ${\cal L}(\omega)=0$, where we set

$${\cal L}(\omega):=||d\theta_\omega||^2_\omega = \int\limits_X d\theta_\omega\wedge\star(d\bar{\theta}_\omega).$$  

\noindent We also know from $(i)$ of Lemma \ref{Lem:general-obstruction} that $d\theta_\omega$ is $\omega$-primitive, so we get

$$0=\Lambda_\omega(d\theta_\omega) = \Lambda_\omega(\partial\theta_\omega^{1,\,0}) + \Lambda_\omega(\partial\theta_\omega^{0,\,1} + \bar\partial\theta_\omega^{1,\,0}) + \Lambda_\omega(\bar\partial\theta_\omega^{0,\,1}) = \Lambda_\omega(\partial\theta_\omega^{0,\,1} + \bar\partial\theta_\omega^{1,\,0}),$$

\noindent where the last identity follows from the previous one for bidegree reasons. We infer that the $(1,\,1)$-form $\partial\theta_\omega^{0,\,1} + \bar\partial\theta_\omega^{1,\,0}$ is $\omega$-primitive. But so are $\partial\theta_\omega^{1,\,0}$ and $\bar\partial\theta_\omega^{0,\,1}$ for bidegree reasons, so we can apply the standard formula (\ref{eqn:prim-form-star-formula-gen}) to get $\star(d\theta_\omega) = \partial\theta_\omega^{1,\,0} -(\partial\theta_\omega^{0,\,1} + \bar\partial\theta_\omega^{1,\,0}) + \bar\partial\theta_\omega^{0,\,1}.$ We infer that \begin{eqnarray}\nonumber d\theta_\omega\wedge\star(d\bar{\theta}_\omega) & = & [\partial\theta_\omega^{1,\,0} +(\partial\theta_\omega^{0,\,1} + \bar\partial\theta_\omega^{1,\,0}) + \bar\partial\theta_\omega^{0,\,1}] \wedge [\partial\theta_\omega^{1,\,0} -(\partial\theta_\omega^{0,\,1} + \bar\partial\theta_\omega^{1,\,0}) + \bar\partial\theta_\omega^{0,\,1}] \\
\nonumber & = & 2\,\partial\theta_\omega^{1,\,0}\wedge\bar\partial\theta_\omega^{0,\,1} - (\partial\theta_\omega^{0,\,1} + \bar\partial\theta_\omega^{1,\,0})^2\end{eqnarray}

\noindent and finally that

\begin{equation}\label{eqn:curly_L_omega}{\cal L}(\omega) = 2\,L(\omega) - \int\limits_X(\partial\theta_\omega^{0,\,1} + \bar\partial\theta_\omega^{1,\,0})^2.\end{equation}

 On the other hand, the Stokes formula implies the first of the following identities

\begin{eqnarray}\label{eqn:integral_dtheta-dtheta}\nonumber 0 & = & \int\limits_X d\theta_\omega\wedge d\theta_\omega = \int\limits_X [\partial\theta_\omega^{1,\,0} +(\partial\theta_\omega^{0,\,1} + \bar\partial\theta_\omega^{1,\,0}) + \bar\partial\theta_\omega^{0,\,1}] \wedge [\partial\theta_\omega^{1,\,0} + (\partial\theta_\omega^{0,\,1} + \bar\partial\theta_\omega^{1,\,0}) + \bar\partial\theta_\omega^{0,\,1}] \\
 & = & 2\,L(\omega) + \int\limits_X(\partial\theta_\omega^{0,\,1} + \bar\partial\theta_\omega^{1,\,0})^2.\end{eqnarray}

 We conclude from (\ref{eqn:curly_L_omega}) and (\ref{eqn:integral_dtheta-dtheta}) that ${\cal L}(\omega)=0$ if and only if $L(\omega)$. Thus, we have proved that $\omega$ is lcK if and only if $L(\omega)=0$, as claimed. 

The identity $L(\omega)= ||\partial\theta_\omega^{1,\,0}||_\omega^2$ follows at once from the general formula (\ref{eqn:prim-form-star-formula-gen}) applied to the primitive $(2,\,0)$-form $\partial\theta_\omega^{1,\,0}$. Indeed, $\star\partial\theta_\omega^{1,\,0} = \partial\theta_\omega^{1,\,0}$, hence $\partial\theta_\omega^{1,\,0}\wedge\bar\partial\theta_\omega^{0,\,1} = \partial\theta_\omega^{1,\,0}\wedge\star\overline{(\partial\theta_\omega^{1,\,0})} = |\partial\theta_\omega^{1,\,0}|^2_\omega\,dV_\omega$.

\vspace{2ex}

$\bullet$ In the case $n\geq 3$, we know from $(ii)$ of Lemma \ref{Lem:general-obstruction} that $\omega$ is lcK if and only if $(d\omega)_{prim}=0$. 

Now, $(d\omega)_{prim} = (\partial\omega)_{prim} + (\bar\partial\omega)_{prim}$ and the forms $(\partial\omega)_{prim}$ and $(\bar\partial\omega)_{prim}$ are conjugate to each other and of different pure types ($(2,\,1)$, respectively $(1,\,2)$), so the vanishing of $(d\omega)_{prim}$ is equivalent to the vanishing of $(\bar\partial\omega)_{prim}$. 

Meanwhile, the standard formula (\ref{eqn:prim-form-star-formula-gen}) applied to the primitive $(2,\,1)$-form $\overline{(\bar\partial\omega)_{prim}} = (\partial\omega)_{prim}$ spells: $$\star\,\overline{(\bar\partial\omega)_{prim}} = i\,\overline{(\bar\partial\omega)_{prim}}\wedge\omega_{n-3}.$$ This proves the identity $L(\omega)=||(\bar\partial\omega)_{prim}||_\omega^2$.

Putting these pieces of information together, we get the following equivalences: \begin{eqnarray*}\omega \hspace{2ex} \mbox{lcK} \iff (d\omega)_{prim}=0 \iff (\bar\partial\omega)_{prim}=0 \iff L(\omega)=0.\end{eqnarray*} The proof is complete. \hfill $\Box$

\section{First variation of the functional: case of complex surfaces}\label{section:1st-variation_surfaces}

Let $S$ be a compact complex surface. (So, we set $X=S$ when $n=2$.) We will compute the differential of the functional $L:{\cal H}_S \longrightarrow[0,\,+\infty)$ defined on the Hermitian cone of $S$. Let $\omega \in {\cal H}_S$. Then, $T_{\omega}{\cal H}_S = C^{\infty} _{1,\,1} (S,\,\mathbb{R})$, so we will compute the differential $$d_\omega L: C^\infty_{1,\,1}(S,\,\R)\longrightarrow\R$$ by computing the derivative of $L(\omega+ t\gamma)$ w.r.t. $t\in(-\varepsilon,\,\varepsilon)$ at $t=0$ for any given real $(1,\,1)$-form $\gamma$.

\begin{Lem}\label{Lem:1st-variation_surfaces} The differential at $\omega$ of the map ${\cal H}_S\ni\omega\mapsto\theta_\omega^{0,\,1}=\Lambda_\omega(\bar\partial\omega)$ is given by $$(d_\omega\theta_\omega^{0,\,1})(\gamma) = \frac{d}{dt}_{|t=0}\Lambda_{\omega+t\gamma}(\bar\partial\omega + t\,\bar\partial\gamma) = \star(\gamma\wedge\star\bar\partial\omega) + \Lambda_\omega(\bar\partial\gamma),$$ while the differential at $\omega$ of $L$ is given by $$(d_\omega L)(\gamma) = 2\,\mbox{Re}\,\int\limits_S\partial\theta_\omega^{1,\,0}\wedge\bar\partial\bigg(\star(\gamma\wedge\star\bar\partial\omega) + \Lambda_\omega(\bar\partial\gamma)\bigg),$$ for every form $\gamma\in C^\infty_{1,\,1}(S,\,\R)$, where $\star=\star_\omega$ is the Hodge star operator defined by the metric $\omega$.

 %(Since $\gamma$ is not assumed of the form $\partial\bar{u} + \bar\partial u$, we tacitly suppose that the definition of $L$ has been extended from ${\cal G}$ to the set of all Hermitian metrics on $S$.) 

\end{Lem}

Before giving the proof of this lemma, we recall the following result from [DP22] that will be used several times in the sequel.

\begin{Lem}\label{Lem:dp22_3.1} ([DP22], Lemmas 3.5 and 3.3) For any complex manifold $X$ of any dimension $n\geq 2$, for any bidegree $(p,\,q)$ and any $C^\infty$ family $(\alpha_t)_{t\in(-\varepsilon,\,\varepsilon)}$ of forms $\alpha_t\in C^\infty_{p,\,q}(X,\,\mathbb{C})$ with $\varepsilon>0$ so small that $\omega+t\gamma>0$ for all $t\in(-\varepsilon,\,\varepsilon)$, the following formulae hold: \begin{eqnarray*}\label{eqn:1st-variation_trace_pq}\frac{d}{dt}\bigg|_{t=0}\,(\Lambda_{\omega+t\gamma}\alpha_t) =  \Lambda_\omega\bigg(\frac{d\alpha_t}{dt}\bigg|_{t=0}\bigg) - (\gamma\wedge\cdot)^\star_\omega\,\alpha_0 =  \Lambda_\omega\bigg(\frac{d\alpha_t}{dt}\bigg|_{t=0}\bigg) + (-1)^{p+q+1}\,\star_\omega(\gamma\wedge\star_\omega\alpha_0).\end{eqnarray*}

\end{Lem} 

 The former of the above equalities appears as such in Lemma 3.5 of [DP22], while the latter equality follows from the former and from formula (27) of Lemma 3.3 of [DP22] which states that $\star_\omega(\eta\wedge\cdot) = (\overline\eta\wedge\cdot)^\star_\omega\,\star_\omega$ for any $(1,\,1)$-form $\eta$ on $X$. Indeed, in our case, taking $\eta = \gamma$ we get $\bar\eta = \gamma$ since $\gamma$ is real. Moreover, composing with $\star_\omega$ on the right and using the standard equality $\star_\omega\star_\omega = (-1)^{p+q}\,\mbox{Id}$ on $(p,\,q)$-forms, we get $\star_\omega(\gamma\wedge\cdot)\star_\omega = (-1)^{p+q}\,(\gamma\wedge\cdot)^\star_\omega$ on $(p,\,q)$-forms.

\vspace{2ex}

\noindent {\it Proof of Lemma \ref{Lem:1st-variation_surfaces}.} The formula for $(d_\omega\theta_\omega^{0,\,1})(\gamma)$ is an immediate consequence of Lemma \ref{Lem:dp22_3.1} applied with $\alpha_t=\bar\partial\omega +t\,\bar\partial\gamma$ (hence also with $(p,\,q)=(1,\,2)$). We further get: \begin{eqnarray}\label{eqn:L_1st_differential1}\nonumber (d_\omega L)(\gamma) & = & \frac{d}{dt}_{|t=0}L(\omega + t\gamma) = \frac{d}{dt}_{|t=0}\int\limits_S\partial\theta^{1,\,0}_{\omega+t\gamma}\wedge\bar\partial\theta^{0,\,1}_{\omega+t\gamma} \\
\nonumber & = & \int\limits_S\partial\bigg(\star(\gamma\wedge\star\partial\omega) + \Lambda_\omega(\partial\gamma)\bigg)\wedge\bar\partial\theta^{0,\,1}_\omega + \int\limits_S \partial\theta^{1,\,0}_\omega\wedge\bar\partial\bigg(\star(\gamma\wedge\star\bar\partial\omega) + \Lambda_\omega(\bar\partial\gamma)\bigg).\end{eqnarray}

\noindent This is the stated formula for $(d_\omega L)(\gamma)$ since the two terms of the r.h.s. expression are mutually conjugated.  \hfill $\Box$

\vspace{3ex}

We will now simplify the above expression of $(d_\omega L)(\gamma)$ starting with a preliminary observation.

\begin{Lem}\label{Lem:simplifying_L_surfaces} Let $(X,\,\omega)$ be an $n$-dimensional complex Hermitian manifold and let $\star=\star_\omega$ be the Hodge star operator defined by $\omega$.

\vspace{1ex}

\noindent $(i)$\, For every $(0,\,1)$-form $\alpha$ on $X$, we have: $$\star(\alpha\wedge\omega) = i\Lambda_\omega(\alpha\wedge\omega_{n-1}).$$

\noindent Moreover, if $n=2$, then $\star(\alpha\wedge\omega)=i\alpha$ for any $(0,\,1)$-form $\alpha$ on $X$.

\vspace{1ex}

\noindent $(ii)$\, If $n=2$, then $\star(\gamma\wedge\alpha) = i\Lambda_\omega(\gamma\wedge\alpha)$ for any $(1,\,1)$-form $\gamma$ and any $(0,\,1)$-form $\alpha$ on $X$. 

\noindent In particular, $\star\bar\partial\omega = i\theta_\omega^{0,\,1}$ for any Hermitian metric $\omega$ on a complex surface.

\vspace{1ex}

\noindent $(iii)$\, In arbitrary dimension $n$, for any $(1,\,1)$-form $\gamma$ and any $(0,\,1)$-form $\alpha$ on $X$, we have: $$\Lambda_\omega(\gamma\wedge\alpha) = (\Lambda_\omega\gamma)\,\alpha + i\,\xi_\alpha\lrcorner\gamma,$$ where $\xi_\alpha$ is the (unique) vector field of type $(1,\,0)$ defined by the requirement $$\xi_\alpha\lrcorner\omega = i\alpha.$$

\end{Lem}

\noindent {\it Proof.} $(i)$\, From the standard formula $\star\Lambda_\omega = L_\omega\star$ (cf. e.g. [Dem97, VI, $\S.5.1$]) we get

\vspace{1ex}

$\Lambda_\omega = \star L_\omega\star$ on even-degreed forms and $\Lambda_\omega = -\star L_\omega\star$ on odd-degreed forms.

\vspace{1ex}

\noindent Consequently, $\star(\alpha\wedge\omega) = \star L_\omega\alpha = -(\star L_\omega\star)\star\alpha = \Lambda_\omega(\star\alpha) = \Lambda_\omega(-(1/i)\,\alpha\wedge\omega^{n-1}/(n-1)!)$, where we used the fact that $\star\star = -1$ on odd-degreed forms and the standard formula (\ref{eqn:prim-form-star-formula-gen}) applied to the (necessarily primitive) $(0,\,1)$-form $\alpha$.

When $n=2$, we get $\star(\alpha\wedge\omega) = i\Lambda_\omega(\alpha\wedge\omega) = i[\Lambda_\omega,\,L_\omega]\,\alpha = -i(1-2)\,\alpha = i\alpha$ after using the general formula $[L_\omega,\,\Lambda_\omega] = (k-n)$ on $k$-forms on $n$-dimensional complex manifolds.

\vspace{1ex}

$(ii)$\, If $n=2$, the map $\omega\wedge\cdot : \Lambda^1T^\star X\longrightarrow \Lambda^3T^\star X$ is an isomorphism at every point of $X$. Since $\gamma\wedge\alpha$ is a $3$-form, there exists a unique $1$-form $\beta$ (necessarily of type $(0,\,1)$) such that $\gamma\wedge\alpha = \omega\wedge\beta$. Moreover, $\beta = \Lambda_\omega(\gamma\wedge\alpha)$ because $\omega\wedge\Lambda_\omega(\gamma\wedge\alpha) = [L_\omega,\,\Lambda_\omega](\gamma\wedge\alpha) = \gamma\wedge\alpha$. Indeed, $\omega\wedge(\gamma\wedge\alpha) = 0$ for bidegree reasons (here $n=2$) and $[L_\omega,\,\Lambda_\omega] = (k-n)$ on $k$-forms.  

Thus, $\gamma\wedge\alpha = \omega\wedge\Lambda_\omega(\gamma\wedge\alpha)$. So, applying $(i)$ for the second identity below, we get: \begin{eqnarray}\nonumber\star(\gamma\wedge\alpha) =\star\bigg(\omega\wedge\Lambda_\omega(\gamma\wedge\alpha)\bigg) = i\Lambda_\omega(\gamma\wedge\alpha).\end{eqnarray}

\noindent To get the last equality, we used (i) for $n=2$ with $\alpha$ replaced by $\Lambda_\omega(\gamma\wedge\alpha)$. 

\vspace{1ex}

In order to prove the formula for $\star\bar\partial\omega$, recall that $\bar\partial\omega = \omega\wedge\theta_\omega^{0,\,1}$, so we get

$$\star\bar\partial\omega = \star(\omega\wedge\theta_\omega^{0,\,1}) = i\theta_\omega^{0,\,1},$$

\noindent where we used again (i) for $n=2$ with $\alpha$ replaced by $\theta_\omega^{0,\,1}$.

\vspace{1ex}

$(iii)$\, Since the claimed identity is pointwise and involves only zero-th order operators, we fix an arbitrary point $x\in X$ and choose local holomorphic coordinates about $x$ such that at $x$ we have

\vspace{1ex}

\hspace{17ex} $\omega = \sum\limits_{a=1}^n idz_a\wedge d\bar{z}_a  \hspace{3ex} \mbox{and} \hspace{3ex} \gamma= \sum\limits_{j=1}^n \gamma_{j\bar{j}}\,idz_j\wedge d\bar{z}_j.$

\vspace{1ex}

\noindent Then, $\Lambda_\omega = -i\sum\limits_{j=1}^n\frac{\partial}{\partial\bar{z}_j}\lrcorner\frac{\partial}{\partial z_j}\lrcorner\cdot$ at $x$. If we set $\alpha = \sum\limits_{j=1}^n\alpha_j\,d\bar{z}_j$ (at any point), we get $\xi_\alpha = \sum\limits_{j=1}^n\alpha_j\,\frac{\partial}{\partial z_j}$ (at $x$) and the following equalities (at $x$): \begin{eqnarray*}\Lambda_\omega(\gamma\wedge\alpha) & = & -i\sum\limits_{j=1}^n\frac{\partial}{\partial\bar{z}_j}\lrcorner\frac{\partial}{\partial z_j}\lrcorner(\gamma\wedge\alpha) \stackrel{(a)}{=} -i\sum\limits_{j=1}^n\frac{\partial}{\partial\bar{z}_j}\lrcorner\bigg(\bigg(\frac{\partial}{\partial z_j}\lrcorner\gamma\bigg)\wedge\alpha\bigg)\\
  & = & -i\sum\limits_{j=1}^n\bigg(\frac{\partial}{\partial\bar{z}_j}\lrcorner\frac{\partial}{\partial z_j}\lrcorner\gamma\bigg)\wedge\alpha + i\sum\limits_{j=1}^n\bigg(\frac{\partial}{\partial z_j}\lrcorner\gamma\bigg)\wedge\bigg(\frac{\partial}{\partial\bar{z}_j}\lrcorner\alpha\bigg)\\
& \stackrel{(b)}{=} & \bigg(\sum\limits_{j=1}^n\gamma_{j\bar{j}}\bigg)\,\alpha - \sum\limits_{j=1}^n\alpha_j\gamma_{j\bar{j}}\,d\bar{z}_j = (\Lambda_\omega\gamma)\,\alpha + i\xi_\alpha\lrcorner\gamma,\end{eqnarray*} where (a) follows from $(\partial/\partial z_j)\lrcorner\alpha = 0$ for bidegree reasons and (b) follows from $(\partial/\partial z_j)\lrcorner\gamma = i\gamma_{j\bar{j}}\,d\bar{z}_j$ and from $(\partial/\partial\bar{z}_j)\lrcorner\alpha = \alpha_j$.

This proves the desired equality at $x$, hence at any point since $x$ was arbitrary. \hfill $\Box$

\vspace{3ex} 

 We can now derive a simplified form of the first variation of the functional $L$.

 \begin{The}\label{The:1st-variation_surfaces} Let $S$ be a compact complex {\bf surface} on which a Hermitian metric $\omega$ has been fixed.

\vspace{1ex}

$(i)$\, The differential at $\omega\in{\cal H}_S$ of the functional $L:{\cal H}_S\longrightarrow[0,\,+\infty)$ evaluated at any form $\gamma\in C^\infty_{1,\,1}(S,\,\R)$ is given by any of the following three formulae: \begin{eqnarray}\label{eqn:first-var_dim_2_1}\nonumber(d_\omega L)(\gamma) & = & -2\,\mbox{Re}\,\int\limits_S\Lambda_\omega(\gamma)\,\partial\theta_\omega^{1,\,0}\wedge\bar\partial\theta_\omega^{0,\,1} - 2\,\mbox{Re}\,\int\limits_S\partial\theta_\omega^{1,\,0}\wedge\bar\partial\Lambda_\omega(\gamma)\wedge\theta_\omega^{0,\,1} + 2\,\mbox{Re}\,\int\limits_S\partial\theta_\omega^{1,\,0}\wedge\bar\partial\Lambda_\omega(\bar\partial\gamma)\\
 & & -2\,\mbox{Re}\,\int\limits_Si\partial\theta_\omega^{1,\, 0}\wedge\bar\partial(\xi_{\theta_\omega^{0,\, 1}}\lrcorner\gamma)\\  
\nonumber  & = & -2\,\mbox{Re}\,\int\limits_S\Lambda_\omega(\gamma)\,|\partial\theta_\omega^{1,\,0}|_\omega^2\,dV_\omega - 2\,\mbox{Re}\,\int\limits_S\partial\theta_\omega^{1,\,0}\wedge\bar\partial\Lambda_\omega(\gamma)\wedge\theta_\omega^{0,\,1} - 2\,\mbox{Re}\,\,i\langle\langle\partial\bar\partial\theta_\omega^{1,\,0},\,\partial\gamma\rangle\rangle_\omega\\
\label{eqn:first-var_dim_2_2} & & -2\,\mbox{Re}\,\int\limits_Si\partial\theta_\omega^{1,\, 0}\wedge\bar\partial(\xi_{\theta_\omega^{0,\, 1}}\lrcorner\gamma)\\
 \label{eqn:first-var_dim_2_3} & = & -2\,\mbox{Re}\,\int\limits_S\partial\theta_\omega^{1,\,0}\wedge\bar\partial\Lambda_\omega(\gamma\wedge\theta_\omega^{0,\,1}) - 2\,\mbox{Re}\,\,i\langle\langle\partial\bar\partial\theta_\omega^{1,\,0},\,\partial\gamma\rangle\rangle_\omega ,\end{eqnarray} where $\star=\star_\omega$ is the Hodge star operator defined by the metric $\omega$ and $\xi_{\theta_\omega^{0,\, 1}}$ is the vector field of type $(1,\,0)$ defined by the requirement $\xi_{\theta_\omega^{0,\, 1}}\lrcorner\omega = i\theta_\omega^{0,\, 1}$.

\vspace{1ex}

$(ii)$\, In particular, for any given $\omega\in{\cal H}_S$, if we choose $\gamma=\partial\theta_\omega^{0,\,1} + \bar\partial\theta_\omega^{1,\,0}$, we have \begin{eqnarray}\nonumber(d_\omega L)(\gamma) = -2\,\mbox{Re}\,\int\limits_Si\partial\theta_\omega^{1,\, 0}\wedge\bar\partial\bigg(\xi_{\theta_\omega^{0,\, 1}}\lrcorner\gamma\bigg) = -2\,\mbox{Re}\,\int\limits_S\partial\theta_\omega^{1,\,0}\wedge\bar\partial\Lambda_\omega(\gamma\wedge\theta_\omega^{0,\,1}).\end{eqnarray}

\end{The}

\noindent {\it Proof.} $(i)$\, From $(ii)$ and $(iii)$ of Lemma \ref{Lem:simplifying_L_surfaces} applied with $\alpha:=i\theta_\omega^{0,\,1}$, we get $$\star(\gamma\wedge\star\bar\partial\omega) = \star(\gamma\wedge i\theta_\omega^{0,\,1}) = i\,\Lambda_\omega(\gamma\wedge i\theta_\omega^{0,\,1}) = -\Lambda_\omega(\gamma)\,\theta_\omega^{0,\,1} - i\xi_{\theta_\omega^{0,\,1}}\lrcorner\gamma.$$

\noindent Formula (\ref{eqn:first-var_dim_2_1}) follows from this and from Lemma \ref{Lem:1st-variation_surfaces}.  

To get (\ref{eqn:first-var_dim_2_2}), we first notice that $\bar\partial\theta_\omega^{0,\,1} = \star\bar\partial\theta_\omega^{0,\,1}$ by the standard formula (\ref{eqn:prim-form-star-formula-gen}) applied to the (necessarily primitive) $(0,\,2)$-form $\bar\partial\theta_\omega^{0,\,1}$. This accounts for the first term on the r.h.s. of (\ref{eqn:first-var_dim_2_2}). Then, we transform the third term on the right-hand side of (\ref{eqn:first-var_dim_2_1}) as follows: \begin{eqnarray*}2\,\mbox{Re}\,\int\limits_S\partial\theta_\omega^{1,\,0}\wedge\bar\partial\Lambda_\omega(\bar\partial\gamma) & \stackrel{(a)}{=} & -2\,\mbox{Re}\,\int\limits_S\partial\theta_\omega^{1,\,0}\wedge\bar\partial\star L_\omega\star(\bar\partial\gamma) \stackrel{(b)}{=} 2\,\mbox{Re}\,\int\limits_S\bar\partial\partial\theta_\omega^{1,\,0}\wedge\star(\omega\wedge\star(\bar\partial\gamma))\\
    & \stackrel{(c)}{=} & 2\,\mbox{Re}\,i\int\limits_S\bar\partial\partial\theta_\omega^{1,\,0}\wedge\star(\bar\partial\gamma) \stackrel{(d)}{=} 2\,\mbox{Re}\,i\int\limits_S\langle\bar\partial\partial\theta_\omega^{1,\,0},\,\partial\bar\gamma\rangle_\omega\,dV_\omega,\end{eqnarray*}

\noindent where we used the standard identity $\Lambda_\omega = -\star L_\omega\star$ on odd-degreed forms to get (a), Stokes to get (b), part $(i)$ of Lemma \ref{Lem:simplifying_L_surfaces} to get (c), and the definition of $\star$ to get (d). Finally, we recall that $\bar\gamma = \gamma$ since $\gamma$ is real.

Finally, (\ref{eqn:first-var_dim_2_3}) follows from Lemma \ref{Lem:1st-variation_surfaces} after using the equality $\star(\gamma\wedge\star\bar\partial\omega)  = -\Lambda_\omega(\gamma\wedge\theta_\omega^{0,\,1})$ (seen above in the proof of (\ref{eqn:first-var_dim_2_1})) and after transforming the third term in (\ref{eqn:first-var_dim_2_1}) as we did above in the proof of (\ref{eqn:first-var_dim_2_2}).

\vspace{1ex}

$(ii)$\, The stated choice of $\gamma$ means that $\gamma$ is the component $(d\theta_\omega)^{1,\,1}$ of type $(1,\,1)$ of the primitive $2$-form $d\theta_\omega$. (See $(i)$ of Lemma \ref{Lem:general-obstruction} for the primitivity statement.) Since $\Lambda_\omega((d\theta_\omega)^{2,\,0})=0$ and $\Lambda_\omega((d\theta_\omega)^{0,\,2})=0$ for bidegree reasons, we infer that \begin{eqnarray*}\Lambda_\omega(\gamma) = \Lambda_\omega((d\theta_\omega)^{1,\,1}) = \Lambda_\omega(d\theta_\omega) = 0.\end{eqnarray*} Therefore, the first two integrals on the r.h.s. of (\ref{eqn:first-var_dim_2_2}) vanish.

Meanwhile, to handle the third integral on the r.h.s. of (\ref{eqn:first-var_dim_2_2}), we notice that $\partial\bar\gamma = \partial\bar\partial\theta_\omega^{1,\,0}$ and this gives the second equality below: \begin{eqnarray*}2\,\mbox{Re}\,\int\limits_S\partial\theta_\omega^{1,\,0}\wedge\bar\partial\Lambda_\omega(\bar\partial\gamma) = 2\,\mbox{Re}\,i\int\limits_S\langle\bar\partial\partial\theta_\omega^{1,\,0},\,\partial\bar\gamma\rangle_\omega\,dV_\omega = -2\,\mbox{Re}\,i||\bar\partial\partial\theta_\omega^{1,\,0}||^2_\omega = 0,\end{eqnarray*} where the first equality above followed from the proof of (\ref{eqn:first-var_dim_2_2}).

Thus, the r.h.s. of formula (\ref{eqn:first-var_dim_2_2}) for $(d_\omega L)(\gamma)$ reduces to its last integral for this choice of $\gamma$. This proves the first claimed equality. 

For the same reason as above, the latter term on the r.h.s. of formula (\ref{eqn:first-var_dim_2_3}) for $(d_\omega L)(\gamma)$ vanishes. This proves the second claimed equality. \hfill $\Box$

\vspace{3ex}

As a first application of (i) of Theorem \ref{The:1st-variation_surfaces}, we deduce the Euler-Lagrange equation for our functional in dimension $2$. The next result can be compared with formula (2.15) of [Vai90]. 

\begin{Cor}\label{Cor:1st-variation_surfaces_cor_E-L_equation} Let $S$ be a compact complex {\bf surface}. The {\bf Euler-Lagrange equation} for the functional $L$ introduced in (i) of Definition \ref{Def:functional_lcK} is \begin{eqnarray*}\label{eqn:surfaces_E-L_equation}\xi_{\theta_\omega^{1,\,0}}\lrcorner\bar\partial\partial\theta_\omega^{1,\,0} + \xi_{\theta_\omega^{0,\,1}}\lrcorner\partial\bar\partial\theta_\omega^{0,\,1} - i\partial^\star\partial\bar\partial\theta_\omega^{1,\,0} + i\bar\partial^\star\bar\partial\partial\theta_\omega^{0,\,1} = 0.\end{eqnarray*} 

\end{Cor}  

\noindent {\it Proof.} $\bullet$ We will use formula (\ref{eqn:first-var_dim_2_3}). The second term on its r.h.s. reads \begin{eqnarray}\label{eqn:1st-variation_surfaces_cor_E-L_equation_proof1}-2\mbox{Re}\,i\langle\langle\partial\bar\partial\theta_\omega^{1,\,0},\,\partial\gamma\rangle\rangle_\omega = \langle\langle - i\partial^\star\partial\bar\partial\theta_\omega^{1,\,0} + i\bar\partial^\star\bar\partial\partial\theta_\omega^{0,\,1},\,\gamma\rangle\rangle_\omega\end{eqnarray} for every real-valued $(1,\,1)$-form $\gamma$.

  $\bullet$ The integral in the first term on the r.h.s. of (\ref{eqn:first-var_dim_2_3}) reads \begin{eqnarray}\label{eqn:1st-variation_surfaces_cor_E-L_equation_proof2}\nonumber\int\limits_S\partial\theta_\omega^{1,\,0}\wedge\bar\partial\Lambda_\omega(\gamma\wedge\theta_\omega^{0,\,1}) & \stackrel{(a)}{=} & \int\limits_S\partial\theta_\omega^{1,\,0}\wedge\star\bigg(\bar\partial\Lambda_\omega(\gamma\wedge\theta_\omega^{0,\,1})\bigg) = \langle\langle\partial\theta_\omega^{1,\,0},\,\partial\Lambda_\omega(\gamma\wedge\theta_\omega^{1,\,0})\rangle\rangle_\omega \\
    & \stackrel{(b)}{=} & -\langle\langle\star\bar\partial\star\partial\theta_\omega^{1,\,0},\,\Lambda_\omega(\gamma\wedge\theta_\omega^{1,\,0})\rangle\rangle_\omega \stackrel{(c)}{=} -\langle\langle\star\bar\partial\partial\theta_\omega^{1,\,0},\,\Lambda_\omega(\gamma\wedge\theta_\omega^{1,\,0})\rangle\rangle_\omega,\end{eqnarray} where (a) and (c) followed from the standard formula (\ref{eqn:prim-form-star-formula-gen}) applied to a $(0,\,2)$-form, resp. a $(2,\,0)$-form, while (b) followed from the standard identity $\partial^\star = -\star\bar\partial\star$.

  Now, $\bar\partial\partial\theta_\omega^{1,\,0}$ is a $(2,\,1)$-form on the $2$-dimensional complex manifold $S$. Since the pointwise map $L_\omega=\omega\wedge\cdot:\Lambda^{1,\,0}T^\star S\longrightarrow\Lambda^{2,\,1}T^\star S$ is bijective (cf. (A) of $\S$\ref{section:standard}), for any $(2,\,1)$-form $\Gamma$ on $S$ there exists a unique $(1,\,0)$-form $\alpha$ such that $\Gamma=\omega\wedge\alpha$. The standard formula (\ref{eqn:prim-form-star-formula-gen}) applied to $\alpha$ yields $\star\alpha = -i\omega\wedge\alpha = -i\Gamma$. Taking $\star$ in the last equality, we get: \begin{eqnarray}\label{eqn:star_10_surfaces_1}\alpha = i\,\star(\omega\wedge\alpha).\end{eqnarray} On the other hand, \begin{eqnarray}\label{eqn:star_10_surfaces_2}\Lambda_\omega(\omega\wedge\alpha) = [\Lambda_\omega,\,L_\omega]\,\alpha = \alpha,\end{eqnarray} where we applied (\ref{eqn:L-Lambda_commutation}) with $n=2$ and $k=1$ to get the last equality. Putting (\ref{eqn:star_10_surfaces_1}) and (\ref{eqn:star_10_surfaces_2}) together, we conclude that \begin{eqnarray*}\star(\omega\wedge\alpha) = -i\,\Lambda_\omega(\omega\wedge\alpha)\end{eqnarray*} for any $(1,\,0)$-form $\alpha$ on a Hermitian complex surface $(S,\,\omega)$.

  In our case, considering the $(2,\,1)$-form $\Gamma:=\bar\partial\partial\theta_\omega^{1,\,0}$ and the unique $(1,\,0)$-form $\alpha$ such that $\Gamma=\omega\wedge\alpha$, we conclude that $\star\bar\partial\partial\theta_\omega^{1,\,0} = \star(\omega\wedge\alpha) = -i\,\Lambda_\omega(\omega\wedge\alpha)$. Hence, (\ref{eqn:1st-variation_surfaces_cor_E-L_equation_proof2}) becomes: \begin{eqnarray*}\label{eqn:1st-variation_surfaces_cor_E-L_equation_proof3}\nonumber\int\limits_S\partial\theta_\omega^{1,\,0}\wedge\bar\partial\Lambda_\omega(\gamma\wedge\theta_\omega^{0,\,1}) & = & i\,\langle\langle\Lambda_\omega(\omega\wedge\alpha),\,\Lambda_\omega(\gamma\wedge\theta_\omega^{1,\,0})\rangle\rangle_\omega = i\,\langle\langle[L_\omega,\,\Lambda_\omega](\omega\wedge\alpha),\,\gamma\wedge\theta_\omega^{1,\,0}\rangle\rangle_\omega \\
    & \stackrel{(a)}{=} & -\langle\langle\bar\partial\partial\theta_\omega^{1,\,0},\,i\theta_\omega^{1,\,0}\wedge\gamma\rangle\rangle_\omega  \stackrel{(b)}{=} -\langle\langle\xi_{\theta_\omega^{1,\,0}}\lrcorner\bar\partial\partial\theta_\omega^{1,\,0},\,\gamma\rangle\rangle_\omega,\end{eqnarray*} where (a) follows from $\omega\wedge\alpha = \bar\partial\partial\theta_\omega^{1,\,0}$ and from formula (\ref{eqn:L-Lambda_commutation}) applied for $k=3$ and $n=2$, while (b) follows from formula (\ref{eqn:xi_bar_contracting_theta01_pre}).

  Thus, the first term on the r.h.s. of (\ref{eqn:first-var_dim_2_3}) reads: \begin{eqnarray}\label{eqn:1st-variation_surfaces_cor_E-L_equation_proof4}-2\,\mbox{Re}\,\int\limits_S\partial\theta_\omega^{1,\,0}\wedge\bar\partial\Lambda_\omega(\gamma\wedge\theta_\omega^{0,\,1}) = \langle\langle\xi_{\theta_\omega^{1,\,0}}\lrcorner\bar\partial\partial\theta_\omega^{1,\,0} + \xi_{\theta_\omega^{0,\,1}}\lrcorner\partial\bar\partial\theta_\omega^{0,\,1} ,\,\gamma\rangle\rangle_\omega,\end{eqnarray} for every real-valued $(1,\,1)$-form $\gamma$.

  $\bullet$ Formulae (\ref{eqn:1st-variation_surfaces_cor_E-L_equation_proof1}) and (\ref{eqn:1st-variation_surfaces_cor_E-L_equation_proof4}) prove the contention.  \hfill $\Box$

\vspace{3ex}

As another application of (i) of Theorem \ref{The:1st-variation_surfaces}, we will now see that the differential $d_\omega L$ vanishes on all the real $(1,\,1)$-forms $\gamma$ that are {\it $\omega$-anti-primitive} (in the sense that $\gamma$ is $\langle\,\,,\,\,\rangle_\omega$-orthogonal to all the $\omega$-primitive $(1,\,1)$-forms, a condition which is equivalent to $\gamma$ being a function multiple of $\omega$). Since $(d_\omega L)(f\omega)$ computes the variation of $L$ in a conformal class, the following statement also follows with no computations from Proposition \ref{Prop:surfaces_scaling-invariant_functional} which shows that the functional $L$ is conformally invariant when $n=2$. However, we prefer giving a direct proof at this point.

\begin{Cor}\label{Cor:1st-variation_surfaces_cor1} Let $S$ be a compact complex {\bf surface} on which a Hermitian metric $\omega$ has been fixed. For any real-valued $C^\infty$ function $f$ on $X$, we have \begin{eqnarray*}(d_\omega L)(f\omega) = 0.\end{eqnarray*}

  In particular, for any real $(1,\,1)$-form $\gamma$ on $S$ we have \begin{eqnarray*}(d_\omega L)(\gamma) = (d_\omega L)(\gamma_{prim}),\end{eqnarray*} where $\gamma_{prim}$ is the $\omega$-primitive component of $\gamma$ in its Lefschetz decomposition.

\end{Cor}

\noindent {\it Proof.} Applying formula (\ref{eqn:first-var_dim_2_1}) with $\gamma = f\omega$ and using the obvious equalities $\Lambda_\omega(f\omega) = 2f$ (recall that $\mbox{dim}_\C S=2$) and $\xi_{\theta_\omega^{0,\, 1}}\lrcorner(f\omega) = f\,(i\theta_\omega^{0,\, 1})$, we get: \begin{eqnarray}\label{eqn:1st-variation_surfaces_cor1_proof}\nonumber(d_\omega L)(f\omega) & = & -4\,\mbox{Re}\,\int\limits_Sf\,\partial\theta_\omega^{1,\,0}\wedge\bar\partial\theta_\omega^{0,\,1} - 4\,\mbox{Re}\,\int\limits_S\partial\theta_\omega^{1,\,0}\wedge\bar\partial f\wedge\theta_\omega^{0,\,1} \\
 \nonumber & & + 2\,\mbox{Re}\,\int\limits_S\partial\theta_\omega^{1,\,0}\wedge\bar\partial\Lambda_\omega(f\bar\partial\omega + \bar\partial f\wedge\omega) -2\,\mbox{Re}\,\int\limits_Si\partial\theta_\omega^{1,\, 0}\wedge(if\bar\partial\theta_\omega^{0,\, 1} + i\bar\partial f\wedge\theta_\omega^{0,\,1}) \\
  & = & T_1 + T_2 + T_3 + T_4,\end{eqnarray} where $T_1$, $T_2$, $T_3$ and $T_4$ stand for the four terms, listed in order, on the r.h.s.  of the above expression for $(d_\omega L)(f\omega)$.

Computing $T_3$, we get: \begin{eqnarray*}T_3 = 2\,\mbox{Re}\,\int\limits_S\partial\theta_\omega^{1,\,0}\wedge\bar\partial(f\,\theta^{0,\,1}_\omega) + 2\,\mbox{Re}\,\int\limits_S\partial\theta_\omega^{1,\,0}\wedge\bar\partial\bigg([\Lambda_\omega,\,L_\omega](\bar\partial f)\bigg),\end{eqnarray*} where we used the equalities $\Lambda_\omega(\bar\partial\omega) = \theta_\omega^{0,\, 1}$ (see (\ref{eqn:Lee-form_formula})) and $\Lambda_\omega(\bar\partial f) = 0$ (which leads to $\Lambda_\omega(\bar\partial f\wedge\omega) = [\Lambda_\omega,\,L_\omega](\bar\partial f)$). Now, it is standard that $[\Lambda_\omega,\,L_\omega] = (n-k)\,\mbox{Id}$ on $k$-forms on an $n$-dimensional complex manifold, so in our case we get $ [\Lambda_\omega,\,L_\omega](\bar\partial f) = \bar\partial f$ since $n=2$ and $k=1$. We conclude that $\bar\partial([\Lambda_\omega,\,L_\omega](\bar\partial f)) = \bar\partial^2f = 0$, hence \begin{eqnarray*}T_3 = 2\,\mbox{Re}\,\int\limits_Sf\,\partial\theta_\omega^{1,\,0}\wedge\bar\partial\theta^{0,\,1}_\omega + 2\,\mbox{Re}\,\int\limits_S\partial\theta_\omega^{1,\,0}\wedge\bar\partial f\wedge\theta^{0,\,1}_\omega =T_4,\end{eqnarray*} where the last equality follows at once from the definition of $T_4$.

Thus, formula (\ref{eqn:1st-variation_surfaces_cor1_proof}) translates to \begin{eqnarray*}(d_\omega L)(f\omega) & = & T_1 + T_2 + T_3 + T_4 \\
  & = & (-4 + 4)\,\mbox{Re}\,\int\limits_Sf\,\partial\theta_\omega^{1,\,0}\wedge\bar\partial\theta_\omega^{0,\,1} + (-4+4)\,\mbox{Re}\,\int\limits_S\partial\theta_\omega^{1,\,0}\wedge\bar\partial f\wedge\theta_\omega^{0,\,1} \\
  & = & 0.\end{eqnarray*} This proves the first statement.

\vspace{1ex}

The second statement follows at once from the first, from the linearity of the map $d_\omega L$ and from the Lefschetz decomposition $\gamma = \gamma_{prim} + (1/2)\,\Lambda_\omega(\gamma)\,\omega$.    \hfill $\Box$

\vspace{2ex}

We hope that it will be possible in the future to prove that any Hermitian metric $\omega$ on a compact complex surface that is a critical point for the functional $L$ is actually an lcK metric.

\section{First variation of the functional: case of dimension $\geq 3$}\label{section:1st-variation_greater than 2} In this section, we suppose that the complex dimension of $X$ is $n\geq 3 $. The goal is to compute the differential of the energy functional $L$ introduced in Definition \ref{Def:functional_lcK}-$(ii)$. Let $\omega$ be a Hermitian metric on $X$ and let $\gamma$ be a real $(1,\,1)$-form. The latter can bee seen as a tangent vector to ${\cal H}_X$ at $\omega$.

\begin{The}\label{The:derivative_L_deg>2} For any Hermitian metric $\omega$ and any real $(1,\,1)$-form $\gamma$, we have: \begin{eqnarray}\label{eqn:derivative_L_deg>2}\nonumber(d_{\omega}L)(\gamma) & = & \int_{X} i(\bar\partial\omega)_{prim}\wedge\overline{(\bar\partial\omega)_{prim}}\wedge\gamma \wedge \omega_{n-4} \\
    & & + 2\text{Re}\,\langle\langle  (\bar\partial\omega)_{prim},\, (\bar\partial\gamma)_{prim}\rangle\rangle_{\omega} - 2\text{Re}\, \langle\langle\theta_\omega^{0,\,1}\wedge\gamma,\,(\bar\partial\omega)_{prim}\rangle\rangle_{\omega}.\end{eqnarray}

\end{The}

\noindent{\it Proof.} Recall (cf. the conjugate of (\ref{eqn:Lee-form_formula_1-0_higher})) that $(n-1)\,\theta_\omega^{0,\,1} = \Lambda_\omega(\bar\partial\omega)$ for any Hermitian metric $\omega$. Now, for any real $t$ sufficiency close to $0$, $\omega+t\gamma$ is again a Hermitian metric on $X$. Taking $\alpha_{t}= \bar\partial\omega+t\,\bar\partial\gamma$ in Lemma \ref{Lem:dp22_3.1}, we get the second equality below: \begin{eqnarray}\label{eqn:derivatie_L_deg>2_proof_1}(n-1)\frac{d}{dt}\bigg|_{t=0}\, \theta^{0,\,1} _{\omega+t\gamma}=\frac{d}{dt}\bigg|_{t=0}\,\Lambda_{\omega+t\gamma}(\bar\partial\omega+t \bar\partial\gamma) = \Lambda_\omega(\bar\partial\gamma)- (\gamma\wedge\cdot)^\star_\omega\,(\bar\partial\omega).\end{eqnarray}

 On the other hand, taking $(d/dt)_{|t=0}$ in the expression for $L(\omega+t\gamma)$ given in (ii) of Definition \ref{Def:functional_lcK} (with $\omega+t\gamma$ in place of $\omega$), we get: \begin{eqnarray}\label{eqn:1st-variation_t_n_geq_3}(d_{\omega}L)(\gamma)=\frac{d}{dt}\bigg|_{t=0}L(\omega+t\gamma)= \frac{d}{dt}\bigg|_{t=0}\int_Xi(\bar\partial\omega+t\bar\partial\gamma)_{prim}\wedge\overline{(\bar\partial\omega+t\bar\partial\gamma)_{prim}}\wedge(\omega+t\gamma)_{n-3},\end{eqnarray} where the subscript {\it prim} indicates the $(\omega + t\gamma)$-primitive part of the form to which it is attached.  

 Now, consider the Lefschetz decompositions (cf. (\ref{eqn:Lefschetz-decomp})) of $\bar\partial\omega$ and $\bar\partial\gamma$ with respect to $\omega$: \begin{eqnarray*}\bar\partial\omega &=& (\bar\partial\omega)_{prim} + \theta_\omega^{0,\,1} \wedge \omega \\
 \bar\partial\gamma &=& (\bar\partial\gamma)_{prim} + \theta_\gamma ^{0,\,1} \wedge \omega\end{eqnarray*} and the Lefschetz decomposition of $\bar\partial\omega+t\bar\partial\gamma$ with respect to $\omega+t\gamma$:
\begin{eqnarray*}\bar\partial\omega+t\bar\partial\gamma &=& (\bar\partial\omega+t\bar\partial\gamma)_{prim} + \theta_{\omega+t\gamma} ^{0,\,1}\wedge(\omega +t \gamma).\end{eqnarray*} 
By the above equations we get:\begin{equation}\label{omega+tgamma prim formula}(\bar\partial\omega+t\bar\partial\gamma)_{prim}= (\bar\partial\omega)_{prim} + \theta_\omega ^{0,\,1} \wedge \omega + t\,(\bar\partial\gamma)_{prim} + t\,\theta_\gamma ^{0,\,1} \wedge \omega - \theta_{\omega+t\gamma} ^{0,\,1} \wedge (\omega +t \gamma),\end{equation} where primitivity is construed w.r.t. the metric $\omega+t\gamma$ in the case of the left-hand side term and w.r.t. the metric $\omega$ in the case of $(\bar\partial\omega)_{prim}$ and $(\bar\partial\gamma)_{prim}$.

%Meanwhile, since $(\bar\partial(\omega+t\gamma))_{prim}$ is $(\omega+t\gamma)$-primitive, by the standard formula (\ref{eqn:prim-form-star-formula-gen}) we get: \begin{eqnarray*}\star \overline{(\bar\partial(\omega+t\gamma))_{prim}}= i \overline{(\bar\partial(\omega+t\gamma))_{prim}}\wedge (\omega+t\gamma)_{n-3},\end{eqnarray*} where $\star$ is taken with respect to $\omega+t\gamma$. From this and from equation (\ref{omega+tgamma prim formula}), we get: \begin{equation*}\label{mega+tgamma prim formula}\star\overline{(\bar\partial(\omega+t\gamma))_{prim}}=i \bigg(\overline{(\bar\partial\omega)_{prim}}+ \overline{\theta_\omega ^{0,\,1}} \wedge \omega + t\overline{(\bar\partial\gamma)_{prim}} + t\overline{\theta_\gamma ^{0,\,1}} \wedge \omega - \overline{\theta_{\omega+t\gamma} ^{0,\,1}} \wedge (\omega +t \gamma)\bigg)\wedge (\omega+t\gamma)_{n-3}.\end{equation*}

Thanks to (\ref{omega+tgamma prim formula}), equality (\ref{eqn:1st-variation_t_n_geq_3}) becomes: 
\begin{eqnarray*}(d_{\omega}L)(\gamma) & = & \frac{d}{dt}_{\bigg|t=0}\int_X i\bigg((\bar\partial\omega)_{prim} + \theta_\omega^{0,\,1} \wedge\omega + t\,(\bar\partial\gamma)_{prim} + t\,\theta_\gamma ^{0,\,1} \wedge \omega - \theta_{\omega+t\gamma} ^{0,\,1} \wedge (\omega +t \gamma)\bigg) \\
\nonumber & &\wedge\bigg(\overline{(\bar\partial\omega)_{prim}}+ \overline{\theta_\omega ^{0,\,1}} \wedge \omega + t\,\overline{(\bar\partial\gamma)_{prim}} + t\,\overline{\theta_\gamma ^{0,\,1}} \wedge \omega - \overline{\theta_{\omega+t\gamma} ^{0,\,1}} \wedge (\omega +t \gamma)\bigg)\wedge (\omega+t\gamma)_{n-3}.
\end{eqnarray*}

Now, \begin{eqnarray*}\frac{d}{dt}_{\bigg|t=0}\bigg(\theta_{\omega+t\gamma}^{0,\,1}\wedge(\omega +t \gamma)\bigg) & = & \theta_\omega^{0,\,1}\wedge\gamma + \bigg(\frac{d}{dt}_{\bigg|t=0}\theta_{\omega+t\gamma}^{0,\,1}\bigg)\wedge\omega \\
  & = & \theta_\omega^{0,\,1}\wedge\gamma + \frac{1}{n-1}\,\bigg(\Lambda_\omega(\bar\partial\gamma) - (\gamma\wedge\cdot)^\star_\omega(\bar\partial\omega)\bigg)\wedge\omega,\end{eqnarray*} where formula (\ref{eqn:derivatie_L_deg>2_proof_1}) was used to get the last equality. Using this, straightforward computations yield: \begin{eqnarray}\label{eqn:first-variaton_L _dim>2_formula}(d_{\omega}L)(\gamma) = I_1 + \overline{I_1} + I_2,\end{eqnarray} where \begin{eqnarray}\label{eqn:first-variaton_L _dim>2_formula_I_2}\nonumber I_2 & = & \int_X i\bigg((\bar\partial\omega)_{prim} + \theta_\omega ^{0,\,1} \wedge\omega - \theta_\omega ^{0,\,1}\wedge\omega\bigg)\wedge\bigg(\overline{(\bar\partial\omega)_{prim}} + \overline{\theta_\omega ^{0,\,1}}\wedge\omega - \overline{\theta_\omega ^{0,\,1}}\wedge\omega\bigg)\wedge\omega_{n-4}\wedge\gamma\\
  & = &  \int_X i(\bar\partial\omega)_{prim}\wedge\overline{(\bar\partial\omega)_{prim}}\wedge\omega_{n-4}\wedge\gamma\end{eqnarray} and \begin{eqnarray}\label{eqn:first-variaton_L _dim>2_formula_I_1}\nonumber I_1 & = & \int_X i\bigg[(\bar\partial\gamma)_{prim} + \theta_\gamma^{0,\,1}\wedge\omega - \theta_\omega^{0,\,1}\wedge\gamma - \frac{1}{n-1}\,\bigg(\Lambda_\omega(\bar\partial\gamma) - (\gamma\wedge\cdot)^\star_\omega(\bar\partial\omega)\bigg)\wedge\omega\bigg]\wedge(\partial\omega)_{prim}\wedge\omega_{n-3}\\
    & = & \int_X i(\bar\partial\gamma)_{prim}\wedge(\partial\omega)_{prim}\wedge\omega_{n-3} - \int_X i\,\theta_\omega^{0,\,1}\wedge\gamma\wedge(\partial\omega)_{prim}\wedge\omega_{n-3},\end{eqnarray} where the last equality follows from $(\partial\omega)_{prim}\wedge\omega_{n-2} = 0$ (a consequence of the $\omega$-primitivity of the $3$-form $(\partial\omega)_{prim}$) which leads to the vanishing of the products of the second and the fourth terms (that are multiples of $\omega$) inside the large parenthesis with $(\partial\omega)_{prim}\wedge\omega_{n-3}$ in the integral on the first line of (\ref{eqn:first-variaton_L _dim>2_formula_I_1}).

  Now, due to the $\omega$-primitivity of the $3$-form $(\partial\omega)_{prim}$, the standard formula (\ref{eqn:prim-form-star-formula-gen}) yields: \begin{eqnarray}\label{eqn:standard-formula_del-omega_app}\star(\partial\omega)_{prim} = i\,(\partial\omega)_{prim}\wedge\omega_{n-3},\end{eqnarray} where $\star = \star_\omega$ is the Hodge star operator induced by $\omega$. Thus, (\ref{eqn:first-variaton_L _dim>2_formula_I_1}) translates to  \begin{eqnarray*}I_1 & = & \int_X (\bar\partial\gamma)_{prim}\wedge\star\overline{(\bar\partial\omega)_{prim}} - \int_X\theta_\omega^{0,\,1}\wedge\gamma\wedge\star\overline{(\bar\partial\omega)_{prim}} \\
    & = & \langle\langle(\bar\partial\gamma)_{prim},\, (\bar\partial\omega)_{prim}\rangle\rangle_{\omega} - \langle\langle\theta_\omega^{0,\,1}\wedge\gamma,\,(\bar\partial\omega)_{prim}\rangle\rangle_{\omega}.\end{eqnarray*}

  This last formula for $I_1$, together with (\ref{eqn:first-variaton_L _dim>2_formula}) and (\ref{eqn:first-variaton_L _dim>2_formula_I_2}), proves the contention.  \hfill $\Box$

  \vspace{3ex}

  The first application of Theorem \ref{The:derivative_L_deg>2} that we give is the computation of the Euler-Lagrange equation for our energy functional $L$ in dimension $n>2$.

  \begin{Cor}\label{Cor:E-L_equation_dim>2} Let $X$ be a compact complex manifold with $\mbox{dim}_\C X=n\geq 3$. The {\bf Euler-Lagrange equation} for the functional $L$ introduced in (ii) of Definition \ref{Def:functional_lcK} is \begin{eqnarray*}\label{eqn:E-L_equation_dim>2}\star\bigg(i(\bar\partial\omega)_{prim}\wedge\overline{(\bar\partial\omega)_{prim}}\wedge\omega_{n-4}\bigg) + \bigg(\bar\partial^\star + i\xi_{\theta^{0,\,1}_\omega}\lrcorner\cdot\bigg)(\bar\partial\omega)_{prim}  + \bigg(\partial^\star - i\xi_{\theta^{1,\,0}_\omega}\lrcorner\cdot\bigg)(\partial\omega)_{prim} = 0,\end{eqnarray*} where the Hodge star operator $\star$, the adjoints and the primitive parts are computed w.r.t. the Hermitian metric $\omega$, the unknown of the equation.  

\end{Cor}

  \noindent {\it Proof.} Using the general formula $\alpha\wedge\beta = \star\alpha\wedge\star\beta$ given in Lemma 5.1. of [Pop22] for any differential forms such that $\mbox{deg}\,\alpha + \mbox{deg}\,\beta = 2n$, the first term on the r.h.s. of (\ref{eqn:derivative_L_deg>2}) transforms as \begin{eqnarray}\label{eqn:E-L_equation_dim>2_proof1}\nonumber\int_X i(\bar\partial\omega)_{prim}\wedge\overline{(\bar\partial\omega)_{prim}}\wedge\gamma \wedge \omega_{n-4} & = & \int_X \star\bigg(i(\bar\partial\omega)_{prim}\wedge\overline{(\bar\partial\omega)_{prim}}\wedge\omega_{n-4}\bigg)\wedge\star\overline\gamma \\
    & = & \bigg\langle\bigg\langle\star\bigg(i(\bar\partial\omega)_{prim}\wedge\overline{(\bar\partial\omega)_{prim}}\wedge\omega_{n-4}\bigg),\,\gamma\bigg\rangle\bigg\rangle\end{eqnarray} for any real $(1,\,1)$-form $\gamma$.

    The second term on the r.h.s. of (\ref{eqn:derivative_L_deg>2}) transforms as \begin{eqnarray}\label{eqn:E-L_equation_dim>2_proof2}\nonumber 2\,\text{Re}\,\langle\langle(\bar\partial\omega)_{prim},\, (\bar\partial\gamma)_{prim}\rangle\rangle_{\omega} & = & 2\text{Re}\,\langle\langle(\bar\partial\omega)_{prim},\,\bar\partial\gamma\rangle\rangle_{\omega} = 2\,\text{Re}\,\langle\langle\bar\partial^\star(\bar\partial\omega)_{prim},\,\gamma\rangle\rangle_{\omega} \\
      & = & \langle\langle\bar\partial^\star(\bar\partial\omega)_{prim} + \partial^\star(\partial\omega)_{prim},\,\gamma\rangle\rangle_{\omega}\end{eqnarray} for any real $(1,\,1)$-form $\gamma$. The first equality above followed from the Lefschetz decomposition $\bar\partial\gamma = (\bar\partial\gamma)_{prim} + \omega\wedge u$ (with some $(0,\,1)$-form $u$) and from $\Lambda_\omega((\bar\partial\omega)_{prim}) =0$.

    The third term on the r.h.s. of (\ref{eqn:derivative_L_deg>2}) transforms, for any real $(1,\,1)$-form $\gamma$, as \begin{eqnarray}\label{eqn:E-L_equation_dim>2_proof3}\nonumber 2\,\text{Re}\,\langle\langle\theta_\omega^{0,\,1}\wedge\gamma,\,(\bar\partial\omega)_{prim}\rangle\rangle_{\omega} & = & 2\,\text{Re}\,\langle\langle(\bar\partial\omega)_{prim},\,\theta_\omega^{0,\,1}\wedge\gamma\rangle\rangle_{\omega} = 2\,\text{Re}\,\langle\langle -i\xi_{\theta^{0,\,1}_\omega}\lrcorner(\bar\partial\omega)_{prim},\,\gamma\rangle\rangle_{\omega} \\
      & = & \langle\langle -i\xi_{\theta^{0,\,1}_\omega}\lrcorner(\bar\partial\omega)_{prim} + i\xi_{\theta^{1,\,0}_\omega}\lrcorner(\partial\omega)_{prim},\,\gamma\rangle\rangle_{\omega},\end{eqnarray} where we used formula (\ref{eqn:xi_bar_contracting_theta01_pre_first}) to get the last equality on the first line.

    \vspace{1ex}

  The contention follows from Theorem \ref{The:derivative_L_deg>2} by putting together (\ref{eqn:E-L_equation_dim>2_proof1}), (\ref{eqn:E-L_equation_dim>2_proof2}) and (\ref{eqn:E-L_equation_dim>2_proof3}). \hfill $\Box$

\vspace{3ex}

Recall that we are interested in the set of critical points of $L$. We now notice that a suitable choice of $\gamma$ in the previous result leads to an explicit description of this set. Since equation (\ref{eqn:derivative_L_deg>2}) is valid for all real $(1,\, 1)$-forms $\gamma$, the choice $\gamma = \omega$ is licit, as any other choice. We get the following

\begin{Cor}\label{Cor:derivative_L_gamma=omega} Let $X$ be a compact complex manifold with $\mbox{dim}_\C X=n\geq 3$ and let $L$ be the functional defined in \ref{Def:functional_lcK}-$(ii)$. For any Hermitian metric $\omega$ on $X$, we have: \begin{equation}\label{eqn:derivative_L_case-gamma=omega}(d_{\omega}L)(\omega)= (n-1)\,\Vert(\bar\partial\omega)_{prim}\Vert_{\omega}^2 = (n-1)\,L(\omega).\end{equation}

\end{Cor}

\noindent{\it Proof.} Taking $\gamma=\omega$ in equation (\ref{eqn:derivative_L_deg>2}), we get:
\begin{eqnarray*}(d_{\omega}L)(\omega)& =& \int_{X} i(\bar\partial\omega)_{prim}\wedge \overline{(\bar\partial\omega)_{prim}}\wedge\omega\wedge\omega_{n-4} +2\text{Re}\,\langle\langle(\bar\partial\omega)_{prim},\,(\bar\partial\omega)_{prim}\rangle\rangle_{\omega} \\
  & & - 2\text{Re}\,\langle\langle \overline{\theta_\omega ^{0,\,1}}\wedge\omega,\,\overline{(\bar\partial\omega)_{prim}} \rangle\rangle_{\omega}\\
  &=& (n-3)i\int_{X} (\bar\partial\omega)_{prim}\wedge\overline{(\bar\partial\omega)_{prim}}\wedge \omega_{n-3} +2\,\Vert(\bar\partial\omega)_{prim}\Vert_{\omega}^2 - 2\text{Re}\,\langle\langle \overline{\theta_\omega ^{0,\,1}},\,\Lambda_\omega((\partial\omega)_{prim})\rangle\rangle_{\omega}\\
  &=& (n-1) \Vert(\bar\partial\omega)_{prim}\Vert_{\omega}^2,\end{eqnarray*} where the last equality followed from $\overline{(\bar\partial\omega)_{prim}}\wedge\omega_{n-3} = -i\,\star\overline{(\bar\partial\omega)_{prim}}$ (see (\ref{eqn:standard-formula_del-omega_app})) and from $\Lambda_\omega((\partial\omega)_{prim})) = 0$ (due to any $\omega$-primitive form lying in the kernel of $\Lambda_\omega$).

\hfill $\Box$\\

An immediate consequence of Corollary \ref{Cor:derivative_L_gamma=omega} is the following 

\begin{Prop}\label{Prop:lcK-critical} Let $X$ be a compact complex manifold with $\mbox{dim}_\C X=n\geq 3$ and let $\omega$ be a Hermitian metric on $X$.

  If $\omega$ is a critical point for the functional $L$ defined in \ref{Def:functional_lcK}-$(ii)$, then $\omega$ is lcK.

\end{Prop}

\noindent {\it Proof.} If $\omega$ is a critical point for $L$, then $(d_{\omega}L)(\gamma) = 0$ for any real $(1,\,1)$-form $\gamma$ on $X$. Taking $\gamma=\omega$ and using (\ref{eqn:derivative_L_case-gamma=omega}), we get $(\bar\partial\omega)_{prim} = 0$. By (ii) of Lemma \ref{Lem:general-obstruction}, this is equivalent to $\omega$ being lcK. \hfill $\Box$

\vspace{2ex}

The converse follows trivially from what we already know. Indeed, if $\omega$ is an lcK metric, $L(\omega) = 0$ (by Lemma \ref{Lem:functional_justification}), so $L$ achieves its minimum at $\omega$ since $L\geq 0$. Any minimum is, of course, a critical point.

\section{Normalised energy functionals when $\mbox{dim}_\C X\geq 3$}\label{section:normalised-functionals}

We start with the immediate observation that the functional introduced in (i) of Definition \ref{Def:functional_lcK} in the case of compact complex surfaces is conformally invariant. In particular, it is scaling-invariant, so it does not need normalising. 

\begin{Prop}\label{Prop:surfaces_scaling-invariant_functional} Let $S$ be a compact complex surface. 

(i)\, For any Hermitian metric $\omega$ on $S$ and any $C^\infty$ function $\lambda:S\longrightarrow(0,\,+\infty)$, the following formula holds: \begin{eqnarray*}\theta_{\lambda\omega}^{1,\,0} = \theta_\omega^{1,\,0} + \frac{1}{\lambda}\,\partial\lambda.\end{eqnarray*}

(ii)\, The functional $L:{\cal H}_S\longrightarrow[0,\,+\infty)$, $L(\omega)=\int_X\partial\theta_\omega^{1,\,0}\wedge\bar\partial\theta_\omega^{0,\,1}$, has the property: \begin{eqnarray*}L(\lambda\omega) = L(\omega)\end{eqnarray*} for every $C^\infty$ function $\lambda:S\longrightarrow(0,\,+\infty)$ and every Hermitian metric $\omega$ on $S$.

\end{Prop}    

\noindent {\it Proof.} (i)\, Recall (cf. (\ref{eqn:Lee-form_formula_1-0})) that $\theta_\omega^{1,\,0} = \Lambda_\omega(\partial\omega)$ and $\theta_\omega^{0,\,1} = \Lambda_\omega(\bar\partial\omega)$. 

On the other hand, for any function $\lambda>0$ and any form $\alpha$ of any bidegree $(p,\,q)$, we have: \begin{eqnarray*}\Lambda_{\lambda\omega}\alpha = \frac{1}{\lambda}\,\Lambda_\omega\alpha,\end{eqnarray*} as can be checked right away. Therefore, we get: \begin{eqnarray*}\theta_{\lambda\omega}^{1,\,0} = \Lambda_{\lambda\omega}\bigg(\partial(\lambda\omega)\bigg) = \frac{1}{\lambda}\,\Lambda_\omega(\lambda\,\partial\omega) + \frac{1}{\lambda}\,\Lambda_\omega(\partial\lambda\wedge\omega) = \theta_\omega^{1,\,0} + \frac{1}{\lambda}\,[\Lambda_\omega,\,L_\omega](\partial\lambda) = \theta_\omega^{1,\,0} + \frac{1}{\lambda}\,\partial\lambda,\end{eqnarray*} where we applied (\ref{eqn:L-Lambda_commutation}) with $n=2$ and $k=1$ to get the last equality. 

\vspace{1ex}

(ii)\, Taking $\partial$ in the formula proved under (i), we get: \begin{eqnarray*}\partial\theta_{\lambda\omega}^{1,\,0} = \partial\theta_\omega^{1,\,0} - \frac{1}{\lambda^2}\,\partial\lambda\wedge\partial\lambda = \partial\theta_\omega^{1,\,0}.\end{eqnarray*} By conjugation, we also get $\bar\partial\theta_{\lambda\omega}^{0,\,1} = \bar\partial\theta_\omega^{0,\,1}$ and the contention follows.  \hfill $\Box$

\vspace{2ex}

By contrast, the functional $L:{\cal H}_X\longrightarrow[0,\,+\infty)$ introduced in (ii) of Definition \ref{Def:functional_lcK} in the case of compact complex manifolds $X$ with $\mbox{dim}_\C X=n\geq 3$ is not scaling-invariant. Indeed, it follows at once from its definition that \begin{eqnarray}\label{eqn:L_homogeneity}L(\lambda\omega) = \lambda^{n-1}\,L(\omega)\end{eqnarray} for every constant $\lambda>0$ and every Hermitian metric $\omega$ on $X$.

  This homogeneity property of $L$ can be used to derive a short proof of the main property of $L$ that was deduced in $\S.$\ref{section:1st-variation_greater than 2} from the result of the computation of the first variation of $L$, namely from Theorem \ref{The:derivative_L_deg>2}.

\begin{Prop}\label{Prop:lcK-critical_bis} (Proposition \ref{Prop:lcK-critical} revisited) Let $X$ be a compact complex manifold with $\mbox{dim}_\C X=n\geq 3$ and let $\omega$ be a Hermitian metric on $X$. The following equivalence holds: \\

\hspace{3ex}  $\omega$ is a {\bf critical point} for the functional $L$ defined in \ref{Def:functional_lcK}-$(ii)$ if and only if $\omega$ is {\bf lcK}.

\end{Prop}

\noindent {\it Proof.} Suppose $\omega$ is a critical point for $L$. This means that $(d_{\omega}L)(\gamma) = 0$ for every real $(1,\,1)$-form $\gamma$ on $X$. Taking $\gamma=\omega$, we get the first eqsuality below: \begin{eqnarray*}0 = (d_{\omega}L)(\omega) = \frac{d}{dt}_{\bigg|t=0}L(\omega + t\omega) = \frac{d}{dt}_{\bigg|t=0}\bigg((1+t)^{n-1}\,L(\omega)\bigg) = (n-1)\,L(\omega).\end{eqnarray*} Thus, whenever $\omega$ is a critical point for $L$, $L(\omega) = 0$. This last fact is equivalent to the metric $\omega$ being lcK thanks to Lemma \ref{Lem:functional_justification}. 

Conversely, if $\omega$ is lcK, it is a minimum point for $L$, hence also a critical point, because $L(\omega) = 0$ by Lemma \ref{Lem:functional_justification}.    \hfill $\Box$

\vspace{2ex}

On the other hand, recall the following by now standard

\begin{Obs}\label{Obs:lcK-bal_one-metric} Let $\omega$ be a Hermitian metric on a complex manifold $X$ with $\mbox{dim}_\C X=n\geq 2$. If $\omega$ is both {\bf lcK} and {\bf balanced}, $\omega$ is {\bf K\"ahler}.

\end{Obs}

\noindent {\it Proof.} The Lefschetz decomposition of $d\omega$ spells $d\omega = (d\omega)_{prim} + \omega\wedge\theta$, where $(d\omega)_{prim}$ is an $\omega$-primitive $3$-form and $\theta$ is a $1$-form on $X$. 

We saw in Lemma \ref{Lem:general-obstruction} that $\omega$ is lcK if and only if $(d\omega)_{prim}=0$. On the other hand, the following equivalences hold: \begin{eqnarray*}\omega \hspace{1ex}\mbox{is balanced}\hspace{1ex}\iff d\omega^{n-1} = 0 \iff \omega^{n-2}\wedge d\omega = 0 \iff d\omega \hspace{1ex}\mbox{is}\hspace{1ex}\omega\mbox{-primitive} \iff d\omega = (d\omega)_{prim}.\end{eqnarray*} 

We infer that, if $\omega$ is both lcK and balanced, $d\omega=0$, so $\omega$ is K\"ahler. \hfill $\Box$

\vspace{2ex}

It is tempting to conjecture the existence of a K\"ahler metric in the more general situation where the lcK and balanced hypotheses are spread over possibly different metrics.

\begin{Conj}\label{Conj:lcK-balanced} Let $X$ be a compact complex manifold with $\mbox{dim}_\C X\geq 3$. If an lcK metric $\omega$ and a balanced metric $\rho$ exist on $X$, there exists a K\"ahler metric on $X$.
  
\end{Conj}

Together with the behaviour of $L$ under rescaling (see (\ref{eqn:L_homogeneity})), this conjecture suggests a natural normalisation for our functional $L$ when $n\geq 3$. 

\begin{Def}\label{Def:normalised-functional_lcK} Let $X$ be a compact complex manifold with $\mbox{dim}_\C X=n\geq 3$. Fix a Hermitian metric $\rho$ on $X$. We define the $\rho$-dependent functional acting on the Hermitian metrics of $X$: \begin{eqnarray}\label{eqn:normalised-functional_lcK}\widetilde{L}_{\rho}: {\cal H}_X\rightarrow[0,\,+\infty), \hspace{5ex} \widetilde{L}_{\rho}(\omega):= \frac{L(\omega)}{\bigg(\int_{X}\omega\wedge\rho_{n-1}\bigg)^{n-1}},\end{eqnarray} where $L$ is the functional introduced in (ii) of Definition \ref{Def:functional_lcK}.

\end{Def}

It follows from (\ref{eqn:L_homogeneity}) that the normalised functional $\widetilde{L}_{\rho}$ is scaling-invariant: \begin{eqnarray*}\widetilde{L}_{\rho}(\lambda\,\omega) = \widetilde{L}_{\rho}(\omega)\end{eqnarray*} for every constant $\lambda>0$. Moreover, thanks to Lemma \ref{Lem:functional_justification}, $\widetilde{L}_{\rho}(\omega)=0$ if and only of $\omega$ is an lcK metric on $X$.

We now derive the formula for the first variation of the normalised functional $\widetilde{L}_{\rho}$ in terms of the similar expression for the unnormalised functional $L$ that was computed in Theorem \ref{The:derivative_L_deg>2}.

\begin{Prop}\label{Prop:1st-variation-normalised-functional_lcK} Let $X$ be a compact complex manifold with $\mbox{dim}_\C X=n\geq 3$. Fix a Hermitian metric $\rho$ on $X$. Then, for any Hermitian metric $\omega$ and any real $(1,\,1)$-form $\gamma$ on $X$, we have:
  \begin{eqnarray}\label{eqn:1st-variation-normalised-functional_lcK}(d_{\omega}\widetilde{L}_{\rho})(\gamma) = \frac{1}{\bigg(\int_{X}\omega\wedge\rho_{n-1}\bigg)^{n-1}}\,\bigg((d_{\omega}L)(\gamma) - (n-1)\,\frac{\int_{X}\gamma\wedge\rho_{n-1}}{\int_{X}\omega\wedge\rho_{n-1}}\,L(\omega)\bigg),\end{eqnarray} where $(d_{\omega}L)(\gamma)$ is given by formula (\ref{eqn:derivative_L_deg>2}) in Theorem \ref{The:derivative_L_deg>2}.

\end{Prop}

\noindent {\it Proof.} Straightforward computations yield: \begin{eqnarray*}(d_{\omega}\widetilde{L}_{\rho})(\gamma) & = & \frac{d}{dt}\bigg[\frac{1}{\bigg(\int_{X}(\omega + t\gamma)\wedge\rho_{n-1}\bigg)^{n-1}}\,L(\omega + t\gamma)\bigg]_{t=0}  = \frac{1}{\bigg(\int_{X}\omega\wedge\rho_{n-1}\bigg)^{n-1}}\,(d_{\omega}L)(\gamma) \\
  & & - \frac{1}{\bigg(\int_{X}\omega\wedge\rho_{n-1}\bigg)^{2(n-1)}}\,(n-1)\,\bigg(\int_{X}\omega\wedge\rho_{n-1}\bigg)^{n-2}\,\bigg(\int_{X}\gamma\wedge\rho_{n-1}\bigg)\,L(\omega).\end{eqnarray*} This is formula (\ref{eqn:1st-variation-normalised-functional_lcK}).  \hfill $\Box$

\vspace{2ex}

A natural question is whether the critical points of any (or some) of the normalised functionals $\widetilde{L}_{\rho}$ are precisely the lcK metrics (if any) on $X$. The following result goes some way in this direction.

\begin{Cor}\label{Cor:omega-critical-for-tilde_gamma-primitive} Let $X$ be a compact complex manifold with $\mbox{dim}_\C X = n\geq 3$. Fix a Hermitian metric $\rho$ on $X$. Suppose a Hermitian metric $\omega$ is a {\bf critical point} for $\widetilde{L}_{\rho}$. Then:

\vspace{1ex}

(i)\, for every {\bf $\rho$-primitive} real $(1,\,1)$-form $\gamma$, $(d_{\omega}L)(\gamma)=0$.

\vspace{1ex}

(ii) if the metric $\rho$ is {\bf Gauduchon}, $(d_{\omega}L)(i\partial\bar\partial\varphi)=0$ for any real-valued $C^2$ function $\varphi$ on $X$. 

\end{Cor}

\noindent {\it Proof.} (i)\, If $\gamma$ is $\rho$-primitive, then $\gamma\wedge\rho_{n-1}=0$, so formula (\ref{eqn:1st-variation-normalised-functional_lcK}) reduces to \begin{eqnarray*}(d_{\omega}\widetilde{L}_{\rho})(\gamma) = \frac{(d_{\omega}L)(\gamma)}{\bigg(\int_{X}\omega\wedge\rho_{n-1}\bigg)^{n-1}}.\end{eqnarray*} Meanwhile, $(d_{\omega}\widetilde{L}_{\rho})(\gamma) = 0$ for every real $(1,\,1)$-form $\gamma$ since $\omega$ is a critical point for $\widetilde{L}_{\rho}$. The contention follows.

\vspace{1ex}

(ii)\, Choose $\gamma: = \omega + i\partial\bar\partial\varphi$ for any function $\varphi$ as in the statement. We get: \begin{eqnarray*}0 \stackrel{(a)}{=} \bigg(\int_{X}\omega\wedge\rho_{n-1}\bigg)^{n-1}\,(d_{\omega}\widetilde{L}_{\rho})(\omega + i\partial\bar\partial\varphi) \stackrel{(b)}{=} (d_{\omega}L)(\omega) - (n-1)\,L(\omega) + (d_{\omega}L)(i\partial\bar\partial\varphi) \stackrel{(c)}{=} (d_{\omega}L)(i\partial\bar\partial\varphi),\end{eqnarray*} where $\omega$ being a critical point for $\widetilde{L}_{\rho}$ gave (a), formula (\ref{eqn:1st-variation-normalised-functional_lcK}) and the metric $\rho$ being Gauduchon (the latter piece of information implying $\int_Xi\partial\bar\partial\varphi\wedge\rho_{n-1} = 0$ thanks to the Stokes theorem) gave (b), while Corollary \ref{Cor:derivative_L_gamma=omega} gave (c).    \hfill $\Box$

\vspace{2ex}

As in the case of surfaces, our hope is that it will be possible in the future to prove that any Hermitian metric $\omega$ on a compact complex manifold of dimension $\geq 3$ that is a critical point for one (or all) of the normalised functionals $\widetilde{L}_\rho$ is actually an lcK metric.

\vspace{2ex}

\noindent {\bf Concluding remarks.} \\

(a)\, Let $X$ be a compact complex manifold with $\mbox{dim}_\C X=n\geq 3$. Fix a Hermitian metric $\rho$ on $X$ and consider the set $U_\rho$ of $\rho$-normalised Hermitian metrics $\omega$ on $X$ such that $$\int\limits_X\omega\wedge\rho_{n-1} = 1.$$ 

By Definition \ref{Def:normalised-functional_lcK}, we have $\widetilde{L}_\rho(\omega) = L(\omega)$ for every $\omega\in U_\rho$. Moreover, since $\widetilde{L}_\rho$ is scaling-invariant, it is completely determined by its restriction to $U_\rho$. Let $$c_\rho:=\inf\limits_{\omega\in{\cal H}_X}\widetilde{L}_\rho(\omega) = \inf\limits_{\omega\in U_\rho}\widetilde{L}_\rho(\omega) = \inf\limits_{\omega\in U_\rho}L(\omega)\geq 0.$$ 

For every $\varepsilon>0$, there exists a Hermitian metric $\omega_\varepsilon\in U_\rho$ such that $c_\rho\leq L(\omega_\varepsilon)< c_\rho + \varepsilon$. Since $U_\rho$ is a relatively compact subset of the space of positive $(1,\,1)$-currents equipped with the weak topology of currents, there exists a subsequence $\varepsilon_k\downarrow 0$ and a positive (see e.g. the terminology of [Dem97, III-1.B.]) $(1,\,1)$-current $T_\rho\geq 0$ on $X$ such that the sequence $(\omega_{\varepsilon_k})_k$ converges weakly to $T_\rho$ as $k\to +\infty$. By construction, we have: $$\int\limits_XT_\rho\wedge\rho_{n-1} = 1.$$ The possible failure of the current $T_\rho\geq 0$ to be either a $C^\infty$ form or strictly positive (for example in the sense that it is bounded below by a positive multiple of a Hermitian metric on $X$) constitutes an obstruction to the existence of minimisers for the functional $\widetilde{L}_\rho$. If it eventually turns out that the critical points of $\widetilde{L}_\rho$, if any, are precisely the lcK metrics of $X$, if any, they will further coincide with the minimisers of $\widetilde{L}_\rho$. In that case, the currents $T_\rho$ will provide obstructions to the existence of lcK metrics on $X$.

\vspace{2ex}

(b)\, The same discussion as in the above (a) can be had on a compact complex surface $S$ using the (already scaling-invariant) functional $L$ introduced in (i) of Definition \ref{Def:functional_lcK} if one can prove that its critical points coincide with the lcK metrics on $S$.

\section*{References}

%\vspace{1ex}

%\noindent [AIOT21]\, D. Angella, N. Istrati, A. Otiman, N. Tardini --- {\it Variational Problems in Conformal Geometry} --- J. Geom. Anal. \textbf{31} (2021), 3230-3251.

\vspace{1ex}

 \noindent [AD15]\, V. Apostolov, G. Dloussky --- {\it Locally Conformally Symplectic Structures on Compact Non-K\"ahler Complex Surfaces} --- Int. Math. Res. Notices, No. {\bf 9} (2016) 2717-2747.

\vspace{1ex}

\noindent [Dem 84]\, J.-P. Demailly --- {\it Sur l'identit\'e de Bochner-Kodaira-Nakano en g\'eom\'etrie hermitienne} --- S\'eminaire d'analyse P. Lelong, P. Dolbeault, H. Skoda (editors) 1983/1984, Lecture Notes in Math., no. {\bf 1198}, Springer Verlag (1986), 88-97.

\vspace{1ex}

\noindent [Dem97]\, J.-P. Demailly --- {\it Complex Analytic and Algebraic Geometry} --- http://www-fourier.ujf-grenoble.fr/~demailly/books.html

\vspace{1ex}

\noindent [DP22]\, S. Dinew, D. Popovici --- {\it A Variational Approach to SKT and Balanced Metrics} --- arXiv:2209.12813v1.

\vspace{1ex}

\noindent [Gau77]\, P. Gauduchon --- {\it Le th\'eor\`eme de l'excentricit\'e nulle} --- C.R. Acad. Sc. Paris, S\'erie A, t. {\bf 285} (1977), 387-390.

\vspace{1ex}

\noindent [Ist19]\. N. Istrati --- {\it Existence Criteria for Special Locally Conformally K\"ahler Metrics} ---  Ann. Mat. Pura Appl. {\bf 198} (2) (2019), 335-353.

\vspace{1ex}

\noindent [Mic83]\, M. L. Michelsohn --- {\it On the Existence of Special Metrics in Complex Geometry} --- Acta Math. {\bf 143} (1983) 261-295.

\vspace{1ex}
 
\noindent [OV22]\, L. Ornea, M. Verbitsky --- {\it Principles of Locally Conformally Kahler Geometry} --- arXiv:2208.07188v2.

\vspace{1ex}

\noindent [Oti14]\, A. Otiman --- {\it Currents on Locally Conformally K\"ahler Manifolds} --- Journal of Geometry and Physics, {\bf 86} (2014), 564-570.

\vspace{1ex}

\noindent [PS22]\, O. Perdu, M. Stanciu --- {\it Vaisman Theorem for lcK Spaces} ---arXiv:2109.01000v3.

\vspace{1ex}

\noindent [Pop22]\, D. Popovici --- {\it Pluriclosed Star Split Hermitian Metrics} --- arXiv e-print DG 2211.10267v1.

\vspace{1ex}

\noindent [Vai76]\, I. Vaisman, --- {\it On Locally Conformal Almost K\"ahler Manifolds} --- Israel J. Math. {\bf 24} (1976) 338-351.

\vspace{1ex}

\noindent [Vai90]\, I. Vaisman, --- {\it On Some Variational Problems for $2$-Dimensional Hermitian Metrics} --- Ann. Global Anal. Geom. {\bf 8}, No. 2 (1990), 137-145.

\vspace{1ex}

\noindent [Voi02]\, C. Voisin --- {\it Hodge Theory and Complex Algebraic Geometry. I.} --- Cambridge Studies in Advanced Mathematics, 76, Cambridge University Press, Cambridge, 2002.

\vspace{6ex}

\noindent Universit\'e Paul Sabatier, Institut de Math\'ematiques de Toulouse

\noindent 118, route de Narbonne, 31062, Toulouse Cedex 9, France

\noindent Email: popovici@math.univ-toulouse.fr \hspace{3ex} and \hspace{3ex} Soheil.Erfan@math.univ-toulouse.fr

\end{document}